\documentclass[11pt,reqno]{amsart}
\usepackage[margin=1.15in]{geometry}
\usepackage{amsmath,amssymb,amsthm,enumitem,mathtools,booktabs}
\usepackage{microtype}
\usepackage{tikz-cd}
\usepackage{hyperref}

\newtheorem{theorem}{Theorem}[section]
\newtheorem{lemma}[theorem]{Lemma}
\newtheorem{proposition}[theorem]{Proposition}
\newtheorem{corollary}[theorem]{Corollary}
\newtheorem{conjecture}[theorem]{Conjecture}

\newtheorem{mainthm}{Theorem}[section]

\theoremstyle{definition}

\theoremstyle{remark}
\newtheorem{remark}[theorem]{Remark}
\numberwithin{equation}{section}

\DeclareMathOperator{\lk}{lk}
\DeclareMathOperator{\Ind}{Ind}
\DeclareMathOperator{\cochord}{cochord}
\DeclareMathOperator{\im}{im}
\DeclareMathOperator{\reg}{reg}

\DeclareMathOperator{\wt}{wt}
\DeclareMathOperator{\Cay}{Cay}
\DeclareMathOperator{\supp}{supp}
\DeclareMathOperator{\depth}{depth}
\DeclareMathOperator{\pd}{pd}
\DeclareMathOperator{\diam}{diam}
\DeclareMathOperator{\rank}{rank}
\newcommand{\X}{\mathcal{X}}
\newcommand{\Y}{\mathcal{Y}}
\newcommand{\coll}{\operatorname{C}}
\newcommand{\ly}{\operatorname{L}}
\newcommand{\VR}{\mathrm{VR}}
\newcommand{\F}{\mathbb{F}_2}
\newcommand{\one}{\mathbf{1}}
\newcommand{\Cl}{\mathrm{Cl}}

\title[Collapsibility of Vietoris--Rips complexes of hypercubes]
      {Collapsibility of Vietoris--Rips complexes of hypercube graphs in codimension two}

\author{S Selvaraja}
\email{selvas@iitbbs.ac.in}
\address{Department of Mathematics,
Indian Institute of Technology Bhubaneswar,
Bhubaneswar-752050, India}


\begin{document}

\begin{abstract}
Shukla conjectured that the collapsibility number of the
Vietoris--Rips complex $\VR(I_n;r)$ of the $n$-dimensional hypercube
graph at scale $r$ equals $2^r$ for all $n\ge r+1$. We prove the conjecture for $n=r+2$ and every $r$, which gives the
first known cases with $r\ge 4$ beyond the cross-polytope case $n=r+1$. Writing $\VR(I_n;r)$ as
the independence complex of the graph $H_{n,r}$ that joins two points
when their Hamming distance exceeds $r$, we construct a cover of
$H_{r+2,r}$ by $2^r$ co-chordal subgraphs. The key general tool is the inequality that, for every graph $G$, the
$\theta$-number of B{\i}y{\i}ko\u{g}lu and Civan of the independence
complex of $G$, and hence the collapsibility number of that complex,
is at most the co-chordal cover number of $G$; this strengthens
Woodroofe's regularity bound. We also show, using a theorem of Alon,
Jin, and Sudakov, that the induced matching number of $H_{n,r}$ equals
$2^r$ for every $n>r$, and that $\VR(I_n;r)$ is vertex decomposable if
and only if $r=1$ or $n=r+1$. As
consequences, we determine the top degree of nonvanishing homology
and the coconnectivity of $\VR(I_{r+2};r)$, and we prove split
injectivity in top-degree integral homology.  Finally, assuming a recent skew Hamming set-pair theorem of Gao and
Zhao, we prove that the collapsibility number of $\VR(I_n;r)$ is at
most $2^{r+1}-4$ for all $n>r\ge 2$.
\end{abstract}
\maketitle

\section{Introduction}
\label{sec:intro}

Let $(X,d)$ be a finite metric space, and let $r\ge 0$. The
\emph{Vietoris--Rips complex} $\VR(X;r)$ is the simplicial complex on
the vertex set $X$ whose faces are the subsets of $X$ of diameter at
most $r$. The $n$-dimensional \emph{hypercube graph} $I_n$ has vertex
set $\F^n$, and two vertices are adjacent if and only if they differ
in exactly one coordinate. The shortest-path metric of $I_n$ is the
Hamming distance
$d_H(\alpha,\beta)=|\{i\in[n]:\alpha_i\ne\beta_i\}|,$
where $[n]=\{1,\ldots,n\}$. We write $\VR(I_n;r)$ for the
Vietoris--Rips complex of $I_n$ with respect to this metric. For
$Y\subseteq\F^n$, we write $\VR(Y;r)$ for the Vietoris--Rips complex
of $Y$ equipped with the restriction of $d_H$.

The topology of these complexes has been studied extensively.
Adamaszek and Adams~\cite[Theorem~1]{AA} proved that $\VR(I_n;2)$ is
homotopy equivalent to a wedge of $3$-spheres, and they asked for the
homotopy type of $\VR(I_n;r)$ for larger values of
$r$~\cite[Question~4]{AA}. For $r=3$, Shukla~\cite[Theorem~A]{Shu}
showed that, for $n\ge 5$, the reduced integral homology of
$\VR(I_n;3)$ is nonzero exactly in degrees $4$ and $7$. Feng~\cite{Fen}
proved further that $\VR(I_n;3)$ is homotopy equivalent to a wedge of
$4$- and $7$-spheres. For general $r$, Adams and
Virk~\cite[Theorem~4.1]{AV} established the lower bound
$\rank\widetilde H_{2^r-1}(\VR(I_n;r);\mathbb Z)
  \ge 2^{n-r-1}\binom{n}{r+1}$ for all integers  $n>r\ge 2.$
Other results on these complexes include connectivity bounds due to
Bendersky and Grbi\'c~\cite{BG}; coconnectivity and cohomological
estimates due to Bendersky, Elia, and Grbi\'c~\cite{BEG}; a study of
facets and rigid facets by Briggs, Feng, and Wells~\cite{BFW}, who
also used Hadamard matrices to establish the nonvanishing of certain
homology groups; and, for $r\le 3$ and $r=n-1$, a decomposition of the
homology over a field of characteristic zero into irreducible
representations of the hyperoctahedral group, due to Galetto,
Monta\~no, and Wellner~\cite{GMW}.

Let $\X$ be a simplicial complex on a finite vertex set $V$. For an
integer $d\ge 0$, an \emph{elementary $d$-collapse} of $\X$ removes a
face $\sigma$ with $|\sigma|\le d$ that is contained in a unique facet
of $\X$, together with all faces of $\X$ that contain $\sigma$; we allow
$\sigma=\emptyset$. The
complex $\X$ is \emph{$d$-collapsible} if it can be reduced to the void
complex by a finite sequence of elementary $d$-collapses. The
\emph{collapsibility number} $\coll(\X)$ is the least $d$ for which
$\X$ is $d$-collapsible. For a field $\Bbbk$, the \emph{Leray number}
$\ly_\Bbbk(\X)$ is the least integer $d$ such that
$\widetilde H_i(\X[W];\Bbbk)=0$ for every $W\subseteq V$ and every
$i\ge d$. Here $\X[W]$ denotes the subcomplex of $\X$ induced by $W$.
Wegner~\cite{Weg} proved that $\ly_\Bbbk(\X)\le\coll(\X)$ for every
field $\Bbbk$, and that a $d$-collapsible complex is homotopy
equivalent to a CW complex of dimension at most $d-1$. Both
invariants play a central role in topological Helly-type
theorems~\cite{Tan}. By Hochster's formula, $\ly_\Bbbk(\X)$ also
equals the Castelnuovo--Mumford regularity of the Stanley--Reisner ring
$S/I_\X$, where $S=\Bbbk[x_v:v\in V]$ and $I_\X$ is the
Stanley--Reisner ideal of $\X$.

Shukla~\cite{Shu} proved that $\coll(\VR(I_n;2))=4$ for $n\ge 3$ and
that $\coll(\VR(I_n;3))=8$ for $n\ge 4$. He conjectured that this
pattern holds for every $r$.

\begin{conjecture}[Shukla~\cite{Shu}]\label{conj:shukla}
For all integers $n$ and $r$ with $n\ge r+1\ge 2$,
\[
  \coll(\VR(I_n;r))=2^r.
\]
\end{conjecture}

The lower bound in Conjecture~\ref{conj:shukla} follows from
coordinate subcubes. Let $Q\subseteq\F^n$ be a coordinate subcube of
dimension $r+1$. Two vertices of $Q$ are at Hamming distance greater
than $r$ if and only if they are antipodal in $Q$. Consequently,
$\VR(Q;r)$ is the boundary of a $2^r$-dimensional cross-polytope.
Since $\VR(Q;r)$ is an induced subcomplex of $\VR(I_n;r)$, we obtain
$2^r\le\ly_\Bbbk(\VR(I_n;r))\le\coll(\VR(I_n;r))$
for every field $\Bbbk$. The difficulty lies in the upper bound. By
Kleitman's diametric theorem~\cite{Kle},
$\dim\VR(I_n;r)=\Theta(n^{\lfloor r/2\rfloor})$ as $n\to\infty$ for
each fixed $r\ge 2$, so the general bound $\coll(\X)\le\dim\X+1$
grows with $n$. Moreover, the coconnectivity results of~\cite{BEG}
concern the homology of the full complex rather than that of its
induced subcomplexes, and so they do not directly bound
collapsibility. Before the present work, Conjecture~\ref{conj:shukla}
was known for $r=1$, for $r\in\{2,3\}$ by the work of
Shukla~\cite{Shu}, and for $n=r+1$, in which case $\VR(I_n;r)$ is the
boundary of a cross-polytope.

Our main result proves Conjecture~\ref{conj:shukla} in codimension
two, which gives the first known cases with $r\ge 4$
beyond the cross-polytope case $n=r+1$. It also determines several related
invariants. Here $\theta$ denotes the $\theta$-number of
B{\i}y{\i}ko\u{g}lu and Civan~\cite{BC}, a recursively defined upper
bound for the collapsibility number, given in terms of deletions and
links (Section~\ref{subsec:invariants}); $d'$ denotes Lew's
invariant~\cite{Lew}; and $H_{n,r}$ denotes the graph on $\F^n$ in
which two vertices are adjacent if and only if their Hamming distance
exceeds $r$.

\begin{mainthm}\label{thmA}
Let $n>r\ge 1$, and let $\Bbbk$ be a field. If $r\in\{1,2\}$ or
$n\in\{r+1,r+2\}$, then
\[
  \reg(S/I_{\VR(I_n;r)})=\ly_\Bbbk(\VR(I_n;r))=\coll(\VR(I_n;r))
  =\theta(\VR(I_n;r))=2^r.
\]
Moreover, $d'(\VR(I_n;r))=2^r$ if $r\le 2$ or $n=r+1$, and
$\cochord(H_{n,r})=2^r$ if $r=1$ or $n\in\{r+1,r+2\}$. In particular,
Conjecture~\ref{conj:shukla} holds for $n=r+2$ and every $r\ge 1$.
\end{mainthm}

Combined with the results of Shukla for $r\le 3$,
Theorem~\ref{thmA} shows that Conjecture~\ref{conj:shukla} remains
open only for $r\ge 4$ and $n\ge r+3$. The smallest open case is
$(n,r)=(7,4)$. For $r=2$, Theorem~\ref{thmA} gives a new proof of
Shukla's equality $\coll(\VR(I_n;2))=4$, and it also determines the
$\theta$-number and Lew's invariant.

Our approach exploits the fact that $\VR(I_n;r)$ is a flag complex:
$\VR(I_n;r)=\Ind(H_{n,r})$, where $\Ind(G)$ denotes the independence
complex of a graph $G$. The \emph{co-chordal cover number}
$\cochord(G)$ is the minimum number of co-chordal subgraphs of $G$
whose edge sets cover $E(G)$. Woodroofe~\cite{W14} proved that
$\ly_\Bbbk(\Ind(G))\le\cochord(G)$ over every field. We strengthen
this bound from the Leray number to the $\theta$-number, and hence to
the collapsibility number.

\begin{mainthm}\label{thmB}
Let $G$ be a graph, and let $\Bbbk$ be a field. Then
\[
  \im(G)\le\ly_\Bbbk(\Ind(G))\le\coll(\Ind(G))\le\theta(\Ind(G))
  \le\min\{d'(\Ind(G)),\cochord(G)\},
\]
where $\im(G)$ denotes the induced matching number of $G$.
Moreover, $d'(\Ind(G))\le\dim\Ind(G)+1$. If $G$ has at least one
edge, then $\Ind(G)$ is homotopy equivalent to a finite CW complex of
dimension at most $\min\{d'(\Ind(G)),\cochord(G)\}-1$.
\end{mainthm}

B{\i}y{\i}ko\u{g}lu and Civan~\cite[Corollary~29]{BC} proved that
$\theta(\X)\le d'(\X)$ for every simplicial complex $\X$; see
also~\cite{Lew,SSS}. The inequality $\theta(\Ind(G))\le\cochord(G)$
is new. It follows from a subadditivity property of the
$\theta$-number under intersections
(Proposition~\ref{prop:theta-subadditivity}), which is the analogue
for $\theta$ of the inequality of Kalai and Meshulam~\cite{KM}
underlying Woodroofe's bound.

Co-chordal covers provide explicit certificates for collapsibility.
For $n=r+2$, we construct a cover of $H_{r+2,r}$ by $2^r$ co-chordal
subgraphs (Theorem~\ref{thm:VR-n-2}); together with
Theorem~\ref{thmB}, this proves Theorem~\ref{thmA} in codimension
two. The cover is optimal, because the lower bound in
Theorem~\ref{thmB} is attained: using the two-sided Hamming set-pair
theorem of Alon, Jin, and Sudakov~\cite{AJS}, we show that
$\im(H_{n,r})=2^r$ for all $n>r\ge 0$
(Corollary~\ref{thm:im-exact}). Thus induced matchings detect exactly
the value conjectured by Shukla for every $n$, and
Theorem~\ref{thmB} shows that Conjecture~\ref{conj:shukla} for a
given pair $(n,r)$ would follow from a cover of $H_{n,r}$ by $2^r$
co-chordal subgraphs. As a combinatorial by-product, we record a
sharp colorful diameter statement for the hypercube, which is a direct
reformulation of the theorem of Alon, Jin, and Sudakov
(Corollary~\ref{cor:colourful-diameter}).

\begin{mainthm}\label{thmC}
Let $n>r\ge 1$.
\begin{enumerate}[label=\textup{(\roman*)}]
\item \textup{(Homotopy dimension.)} If $r\le 2$ or $n\le r+2$, then
for every nonempty $Y\subseteq\F^n$ the complex $\VR(Y;r)$ is
homotopy equivalent to a finite CW complex of dimension at most
$2^r-1$.

\item \textup{(Top homology.)} If $r\le 3$ or $n\le r+2$, then $2^r-1$
is the largest degree in which $\VR(I_n;r)$ has nonzero reduced
integral homology, and the coconnectivity of $\VR(I_n;r)$ equals
$2^r$.

\item \textup{(Spheres from subcubes.)} If $Y\subseteq\F^n$ contains an
$(r+1)$-dimensional coordinate subcube, then
$\widetilde H_{2^r-1}(\VR(Y;r);\mathbb Z)$ contains $\mathbb Z$ as a
direct summand.

\item \textup{(Split injectivity.)} If $r\le 2$ or $n\le r+2$, then
every inclusion $Y'\subseteq Y\subseteq\F^n$ induces a split
injection
\[
  \widetilde H_{2^r-1}(\VR(Y';r);\mathbb Z)\longrightarrow
  \widetilde H_{2^r-1}(\VR(Y;r);\mathbb Z).
\]

\item \textup{(Vertex decomposability.)} The following conditions are
equivalent:
\begin{enumerate}[label=\textup{(\alph*)}]
\item $\VR(I_n;r)$ is vertex decomposable;
\item $\VR(I_n;r)$ is Cohen--Macaulay over some field;
\item $\VR(I_n;r)$ is pure;
\item $r=1$ or $n=r+1$.
\end{enumerate}
Moreover, $2\le\depth(S/I_{\VR(I_n;r)})\le 2^r$ over every field.
\end{enumerate}
\end{mainthm}

Part~(ii) confirms, for $n=r+2$, the part of a conjecture of
Shukla~\cite{Shu} concerning the top degree of nonvanishing homology,
and closes the gap between the upper bounds of~\cite[Theorem~4.2]{BEG}
and the nonvanishing recorded in~\cite[Table~2]{BEG}. The direct
summand in part~(iii) arises from a simplicial retraction onto a
coordinate subcube (Proposition~\ref{prop:retraction}), and the split
injectivity in part~(iv) follows from a bound on the $\theta$-number
of punctured Hamming balls (Lemma~\ref{lem:punctured}) combined with
the Mayer--Vietoris sequence. Part~(v) shows that vertex
decomposability cannot be used to prove Shukla's conjecture in the
remaining cases: for $r\ge 2$ and $n\ge r+2$, the complex
$\VR(I_n;r)$ has no shedding vertices, although
$\theta(\VR(I_{r+2};r))=2^r$ by Theorem~\ref{thmA}.

For $n\ge r+3$, our bounds come from Lew's invariant. By
Lemma~\ref{lem:dprime-VR}, $d'(\VR(I_n;r))$ is the largest integer
$k$ for which there exist pairwise $r$-close points
$v_1,\ldots,v_k\in\F^n$ and points $w_1,\ldots,w_k\in\F^n$ such that
\begin{equation}\label{eq:skew-intro}
  d_H(w_i,v_i)>r ~ (1\le i\le k)
  ~\text{and}~
  d_H(w_i,v_j)\le r ~ (1\le j<i\le k).
\end{equation}
If the requirement that the points $v_i$ be pairwise $r$-close is
dropped, such a configuration is a skew Hamming set-pair system.
Alon, Jin, and Sudakov~\cite{AJS} determined the maximum size of a
symmetric Hamming set-pair system and asked whether the same bound
holds in the skew setting. Gao and Zhao~\cite{GZ} recently answered
this question in the affirmative, obtaining the sharp bound
$2^{r+1}$. As~\cite{GZ} is currently available only as a preprint,
we isolate the statement that we need as a hypothesis, and we
indicate explicitly every result that depends on it.

\medskip
\noindent\textbf{Hypothesis $(\mathrm S_r)$.}
\emph{For every integer $n>r$, if $(a_1,b_1),\ldots,(a_k,b_k)$ are
pairs of points of $\F^n$ such that $d_H(a_i,b_i)>r$ for every $i$
and $d_H(a_i,b_j)\le r$ whenever $i<j$, then $k\le 2^{r+1}$.}

\medskip
\noindent
Hypothesis~$(\mathrm S_r)$ is precisely the binary case, with $t=r$,
of the theorem of Gao and Zhao~\cite[Theorem~2.2]{GZ}, which we
recall as Theorem~\ref{thm:GZ}. Consequently, $(\mathrm S_r)$ is not
an additional conjecture: it holds for every $r\ge 1$ by~\cite{GZ},
and all results stated below under $(\mathrm S_r)$ are unconditional
once that theorem is accepted.

\medskip
\noindent
Let $\kappa(n,r)=\dim\VR(I_n;r)+1$, which is given explicitly by
Kleitman's diametric theorem~\cite{Kle}:
\[
  \kappa(n,r)=
  \begin{cases}
    \displaystyle\sum_{i=0}^{r/2}\binom{n}{i}
      & \text{if $r$ is even},\\[3mm]
    \displaystyle 2\sum_{i=0}^{(r-1)/2}\binom{n-1}{i}
      & \text{if $r$ is odd}.
  \end{cases}
\]
For $n>r\ge 2$, set $U(n,r)=\min\{\kappa(n,r),\,2^{r+1}-4\}$.

\begin{mainthm}\label{thmD}
Let $n>r\ge 2$, and let $\Bbbk$ be a field.
\begin{enumerate}[label=\textup{(\alph*)}]
\item \textup{(Unconditional.)} We have
\begin{equation}\label{eq:best-intro}
  2^r=\im(H_{n,r})\le\ly_\Bbbk(\VR(I_n;r))\le\coll(\VR(I_n;r))
  \le\theta(\VR(I_n;r))\le d'(\VR(I_n;r))\le\kappa(n,r).
\end{equation}
\item \textup{(Conditional on $(\mathrm S_r)$.)} If
Hypothesis~$(\mathrm S_r)$ holds, then
$d'(\VR(I_n;r))\le\min\{\kappa(n,r),2^{r+1}\}$ and
\begin{equation}\label{eq:best-intro-cond}
  \theta(\VR(I_n;r))\le U(n,r)\le 2^{r+1}-4.
\end{equation}
Consequently, every complex $\VR(Y;r)$ with
$\emptyset\ne Y\subseteq\F^n$ is homotopy equivalent to a finite CW
complex of dimension at most $2^{r+1}-5$.
\end{enumerate}
\end{mainthm}

Under $(\mathrm S_r)$, part~(b) determines the value conjectured by
Shukla up to a factor strictly less than $2$, uniformly in $n$. The
bound~\eqref{eq:best-intro-cond} is obtained from Lew's invariant by
a link-padding argument (Theorem~\ref{thm:VR-padded}). It holds
unconditionally when $r=2$ or $n\le n_0(r)$, where $n_0(r)$ is given
in Corollary~\ref{cor:best-bound}. For example, we have
unconditionally $16\le\coll(\VR(I_7;4))\le 29$, and under
$(\mathrm S_4)$ the upper bound improves to $28$. The bound $2^{r+1}$
in $(\mathrm S_r)$ is sharp for general skew set-pair systems, but
the extremal construction does not satisfy the requirement
in~\eqref{eq:skew-intro} that the points $v_i$ be pairwise $r$-close.
We conjecture that this requirement improves the bound to $2^r$
(Conjecture~\ref{conj:skew}), which would imply Shukla's conjecture.

\vskip 1mm \noindent
\emph{Organization:}
Section~\ref{sec:prelim} collects preliminaries on simplicial
complexes, the $\theta$-number, Stanley--Reisner rings, co-chordal
covers, and Lew's invariant, and proves Theorem~\ref{thmB}.
Section~\ref{sec:VR} introduces the Vietoris--Rips complexes of
hypercubes, determines the induced matching number of $H_{n,r}$, and
records the colorful diameter statement. Section~\ref{sec:small} proves
Theorem~\ref{thmA}. Section~\ref{sec:uniform} proves
Theorem~\ref{thmD}, separating the unconditional bounds from those
that depend on Hypothesis~$(\mathrm S_r)$. Finally,
Section~\ref{sec:topology} proves Theorem~\ref{thmC}.

\section{Preliminaries}
\label{sec:prelim}

All graphs are finite and simple. For a graph $G$ and a vertex
$v\in V(G)$, we write $N(v)=N_G(v)$ and $N[v]=N(v)\cup\{v\}$ for the
open and closed neighborhoods of $v$, respectively, and
$\deg_G(v)=|N(v)|$. For $U\subseteq V(G)$, we write $G[U]$ for the
subgraph of $G$ induced by $U$, and $G-U=G[V(G)\setminus U]$. The
independence number of $G$ is denoted by $\alpha(G)$, and the
complement of $G$ by $\overline G$. An \emph{induced matching} in $G$ is a
matching $M$ such that the subgraph induced by the endpoints of the
edges of $M$ has edge set $M$. The \emph{induced matching number}
$\im(G)$ is the maximum size of an induced matching in $G$.

\subsection{Simplicial complexes}
\label{subsec:complexes}

A \emph{simplicial complex} on a finite ground set $V$ is a nonempty
family $\X\subseteq 2^V$ that is closed under taking subsets. The
empty family is called the \emph{void complex} and is also denoted
by $\emptyset$; it is not a simplicial complex in the above sense,
and it arises only as the result of a complete collapse and as the
link of a ghost vertex (see below). The elements of $\X$ are its
\emph{faces}, and its maximal faces are its \emph{facets}. The
\emph{vertices} of $\X$ are the elements $v\in V$ with $\{v\}\in\X$;
the remaining elements of $V$ are \emph{ghost vertices}. We set
$\dim\X=\max_{\sigma\in\X}|\sigma|-1$, and for $W\subseteq V$ we write
$2^W$ for the simplex on $W$. In particular,
$2^\emptyset=\{\emptyset\}$ has dimension $-1$ and is distinct from
the void complex.
All invariants considered in this paper are unchanged by adding or
removing ghost vertices, and we therefore compare complexes as
families of sets.

For $U\subseteq V$, the subcomplex of $\X$ \emph{induced} by $U$ is
$\X[U]=\{\sigma\in\X:\sigma\subseteq U\}$. For a vertex $v$ of $\X$,
the \emph{deletion} and the \emph{link} of $v$ are
\[
  \X\setminus v=\X[V\setminus\{v\}]
  ~\text{and}~
  \lk(\X,v)=\{\sigma\in\X: v\notin\sigma,\ \sigma\cup\{v\}\in\X\},
\]
respectively. The \emph{join} of complexes $\X$ and $\Y$ on disjoint
ground sets is $\X*\Y=\{\sigma\cup\eta:\sigma\in\X,\ \eta\in\Y\}$. A
vertex $v$ is a \emph{cone point} of $\X$ if
$\X=2^{\{v\}}*(\X\setminus v)$. A complex is \emph{flag} if every
clique in its $1$-skeleton is a face; equivalently, every minimal
non-face contained in its vertex set has exactly two elements.

The \emph{clique complex} $\Cl(\Gamma)$ of a graph $\Gamma$ is the
complex of cliques of $\Gamma$, and the \emph{independence complex}
of a graph $G$ is $\Ind(G)=\Cl(\overline G)$. Every flag complex is
the independence complex of a graph on its vertex set, namely the
complement of its $1$-skeleton. For a vertex $v$ of $G$ and for
graphs $G_1$ and $G_2$ on disjoint vertex sets, we have
\begin{equation}\label{eq:ind-link}
\begin{gathered}
  \Ind(G)\setminus v=\Ind(G-v),~
  \lk(\Ind(G),v)=\Ind(G-N[v]),\\
  \Ind(G_1\sqcup G_2)=\Ind(G_1)*\Ind(G_2).
\end{gathered}
\end{equation}

\subsection{Leray number, collapsibility number, and the
\texorpdfstring{$\theta$}{theta}-number}
\label{subsec:invariants}

Fix a field $\Bbbk$. Unless stated otherwise, homology is taken with
coefficients in $\Bbbk$. Let $\X$ be a simplicial complex on $V$.
\begin{enumerate}[label=\textup{(\alph*)}]
\item The \emph{Leray number} $\ly(\X)=\ly_\Bbbk(\X)$ is the least
integer $d\ge 0$ such that $\widetilde H_i(\X[U];\Bbbk)=0$ for every
$U\subseteq V$ and every $i\ge d$. We say that $\X$ is
\emph{$d$-Leray} if $\ly(\X)\le d$.
\item Let $d\ge 0$. An \emph{elementary $d$-collapse} of $\X$ removes
a face $\sigma$ with $|\sigma|\le d$ that is contained in a unique
facet of $\X$, together with all faces of $\X$ that contain $\sigma$;
we allow $\sigma=\emptyset$. The complex $\X$ is
\emph{$d$-collapsible} if it can be reduced to the void complex by a
finite sequence of elementary $d$-collapses. The
\emph{collapsibility number} $\coll(\X)$ is the least $d$ for which
$\X$ is $d$-collapsible. In particular, $\coll(2^W)=0$ for every
$W\subseteq V$.
\item The \emph{$\theta$-number} of B{\i}y{\i}ko\u{g}lu and
Civan~\cite{BC} is defined recursively as follows. If $\X$ is a
simplex, then $\theta(\X)=0$. Otherwise,
\begin{equation}\label{eq:theta-rec}
  \theta(\X)=\min_v\max\{\theta(\X\setminus v),\,
  1+\theta(\lk(\X,v))\},
\end{equation}
where the minimum is taken over all vertices $v$ of $\X$ that are not
cone points.
\end{enumerate}

We shall use the following standard facts.
\begin{enumerate}[label=\textup{(F\arabic*)},leftmargin=3em]
\item\label{F:leray}
$\ly(\X)\le\coll(\X)\le\theta(\X)$~\cite{Weg,BC}.
\item\label{F:mono}
If $U\subseteq V$, then $\ly(\X[U])\le\ly(\X)$,
$\coll(\X[U])\le\coll(\X)$, and $\theta(\X[U])\le\theta(\X)$. The
first inequality is immediate from the definition, the second is due
to Tancer~\cite{Tan}, and the third is~\cite[Proposition~2]{BC}.
\item\label{F:join}
If $\X$ and $\Y$ are complexes on disjoint ground sets, then
$\theta(\X*\Y)=\theta(\X)+\theta(\Y)$~\cite[Theorem~11]{BC}. In
particular, $\theta(\X*2^A)=\theta(\X)$ for every finite set $A$
disjoint from the ground set of $\X$.
\item\label{F:glue}
For every vertex $v$ of $\X$, we have
$\X=(\X\setminus v)\cup\bigl(2^{\{v\}}*\lk(\X,v)\bigr)$, and the two
subcomplexes on the right intersect in $\lk(\X,v)$~\cite{Bj}.
\end{enumerate}

The next two lemmas are essentially contained in~\cite{BC}:
Lemma~\ref{lem:theta-basic} follows from~\cite[Corollary~5]{BC}, and
Lemma~\ref{lem:theta-link-mono} follows from~\cite[Theorem~4]{BC}
together with~\ref{F:mono}. We include short proofs for completeness.

\begin{lemma}\label{lem:theta-basic}
For every vertex $v$ of $\X$,
$\theta(\X)\le\max\{\theta(\X\setminus v),\,1+\theta(\lk(\X,v))\}.$
Moreover, $\theta(\X)\le\dim\X+1$.
\end{lemma}

\begin{proof}
If $v$ is not a cone point, the first inequality follows directly
from~\eqref{eq:theta-rec}. If $v$ is a cone point, then
$\theta(\X)=\theta(\X\setminus v)$ by~\ref{F:join}. The second
assertion is~\cite[Corollary~8]{BC}.
\end{proof}

\begin{lemma}\label{lem:theta-link-mono}
For every vertex $v$ of $\X$, we have
$\theta(\lk(\X,v))\le\theta(\X)$.
\end{lemma}

\begin{proof}
We proceed by induction on the number of vertices of $\X$. If $\X$ is
a simplex, then so is $\lk(\X,v)$, and the claim holds. Otherwise,
let $u$ be a vertex at which the minimum in~\eqref{eq:theta-rec} is
attained, so that $\theta(\X\setminus u)\le\theta(\X)$ and
$1+\theta(\lk(\X,u))\le\theta(\X)$. If $u=v$, the claim follows, so
we may assume that $u\ne v$.
Suppose first that $\{u,v\}\notin\X$. Then
$\lk(\X,v)=\lk(\X\setminus u,v)$, and the induction hypothesis
applied to $\X\setminus u$ gives
$\theta(\lk(\X,v))\le\theta(\X\setminus u)\le\theta(\X)$.
Suppose now that $\{u,v\}\in\X$. Then $u$ is a vertex of
$\lk(\X,v)$, $v$ is a vertex of $\lk(\X,u)$, and
\[
  \lk(\X,v)\setminus u=\lk(\X\setminus u,v),
  ~
  \lk(\lk(\X,v),u)=\lk(\lk(\X,u),v).
\]
Applying Lemma~\ref{lem:theta-basic} to $\lk(\X,v)$ at $u$, and the
induction hypothesis to $\X\setminus u$ and to $\lk(\X,u)$, we obtain
\[
  \theta(\lk(\X,v))
  \le\max\{\theta(\X\setminus u),\,1+\theta(\lk(\X,u))\}
  =\theta(\X).
\]
This completes the proof.
\end{proof}

\begin{lemma}\label{lem:certificate}
Let $d\ge 1$, let $v_1,\ldots,v_m$ be distinct vertices of $\X$, and
set $\X_0=\X$ and $\X_i=\X_{i-1}\setminus v_i$ for $1\le i\le m$. If
$\theta(\lk(\X_{i-1},v_i))\le d-1$ for $1\le i\le m$ and
$\theta(\X_m)\le d$, then $\theta(\X)\le d$.
\end{lemma}

\begin{proof}
Apply Lemma~\ref{lem:theta-basic} to $\X_{i-1}$ at $v_i$ for
$i=m,m-1,\ldots,1$ in turn.
\end{proof}

\begin{lemma}\label{lem:link}
Let $\X$ be a simplicial complex that is not a simplex. Then
\[
  \theta(\X)\le 1+\max_v\theta(\lk(\X,v))
  ~\text{and}~
  \ly(\X)\le 1+\max_v\ly(\lk(\X,v)),
\]
where both maxima are taken over all vertices $v$ of $\X$. If the
automorphism group of $\X$ acts transitively on its vertices, then
each maximum is attained at every vertex.
\end{lemma}

\begin{proof}
Let $d=1+\max_v\theta(\lk(\X,v))$, let $v_1,\ldots,v_m$ be the
vertices of $\X$, and define $\X_i$ as in
Lemma~\ref{lem:certificate}. Each $\lk(\X_{i-1},v_i)$ is an induced
subcomplex of $\lk(\X,v_i)$, so
$\theta(\lk(\X_{i-1},v_i))\le d-1$ by~\ref{F:mono}. Moreover,
$\X_m=\{\emptyset\}$. Hence $\theta(\X)\le d$ by
Lemma~\ref{lem:certificate}.
Now let $d=1+\max_v\ly(\lk(\X,v))$. We show by induction on $|W|$
that $\widetilde H_i(\X[W])=0$ for every $W\subseteq V$ and every
$i\ge d$. If $W$ contains no vertex of $\X$, then
$\X[W]=\{\emptyset\}$ and the claim holds because $d\ge 1$.
Otherwise, let $v\in W$ be a vertex of $\X$. The link of $\X[W]$ at
$v$ is $\lk(\X,v)[W\setminus\{v\}]$. Applying~\ref{F:glue} to
$\X[W]$ and using the Mayer--Vietoris sequence, we obtain the exact
sequence
\[
  \widetilde H_i(\X[W\setminus\{v\}])\longrightarrow
  \widetilde H_i(\X[W])\longrightarrow
  \widetilde H_{i-1}\bigl(\lk(\X,v)[W\setminus\{v\}]\bigr).
\]
For $i\ge d$, the left-hand term vanishes by the induction
hypothesis, and the right-hand term vanishes because
$\ly(\lk(\X,v))\le d-1$. Hence the middle term vanishes.
The final assertion holds because automorphisms of $\X$ map links to
isomorphic links.
\end{proof}

\begin{proposition}\label{prop:theta-subadditivity}
Let $\X_1,\ldots,\X_c$ be simplicial complexes on a common finite
ground set $V$. Then
\[
  \theta\Bigl(\bigcap_{j=1}^c\X_j\Bigr)\le\sum_{j=1}^c\theta(\X_j).
\]
\end{proposition}

\begin{proof}
Set $\X=\bigcap_{j=1}^c\X_j$ and $s=\sum_{j=1}^c\theta(\X_j)$. We
proceed by induction on $|V|$. If $\X$ is a simplex, there is nothing
to prove. If $\X$ has a ghost vertex $v$, then
$\X=\bigcap_{j=1}^c\X_j[V\setminus\{v\}]$, and the induction
hypothesis together with~\ref{F:mono} gives $\theta(\X)\le s$.
We may therefore assume that every element of $V$ is a vertex of
$\X$. If every $\X_j$ were a simplex $2^{W_j}$, then
$\X=2^{\bigcap_j W_j}$ would be a simplex; hence some $\X_i$ is not
a simplex. Let $v$ be a vertex at which the minimum
in~\eqref{eq:theta-rec} is attained for $\X_i$. Then
$\theta(\X_i\setminus v)\le\theta(\X_i)$ and
$\theta(\lk(\X_i,v))\le\theta(\X_i)-1$. Since deletion and link
commute with intersections, we have
$\X\setminus v=\bigcap_{j=1}^c(\X_j\setminus v)$
and
$\lk(\X,v)=\bigcap_{j=1}^c\lk(\X_j,v).$
By the induction hypothesis,~\ref{F:mono}, and
Lemma~\ref{lem:theta-link-mono}, it follows that
$\theta(\X\setminus v)\le s$ and $\theta(\lk(\X,v))\le s-1$.
Lemma~\ref{lem:theta-basic} now yields $\theta(\X)\le s$.
\end{proof}

\begin{remark}
Proposition~\ref{prop:theta-subadditivity} is the analogue for the
$\theta$-number of the inequality $\ly(\X\cap\Y)\le\ly(\X)+\ly(\Y)$
of Kalai and Meshulam~\cite{KM}, which underlies Woodroofe's
regularity bound~\cite{W14}. Combined with~\ref{F:leray}, it gives
$\coll(\bigcap_j\X_j)\le\sum_j\theta(\X_j)$. We do not know whether
the collapsibility number itself is subadditive under intersections.
\end{remark}

The following lemma is a consequence of Wegner's theorem that a
$d$-collapsible complex is homotopy equivalent to a CW complex of
dimension at most $d-1$~\cite{Weg}.

\begin{lemma}\label{lem:theta-homotopy}
Let $d\ge 1$. If $\theta(\X)\le d$, then the geometric realization
$|\X|$ is homotopy equivalent to a finite CW complex of dimension at
most $d-1$. Consequently, $\widetilde H_i(\X;\mathbb Z)=0$ for every
$i\ge d$, and $\widetilde H_{d-1}(\X;\mathbb Z)$ is free abelian.
\end{lemma}

\begin{proof}
By~\ref{F:leray}, $\coll(\X)\le\theta(\X)\le d$, so the first
assertion follows from Wegner's theorem~\cite{Weg}. The cellular
chain complex of a CW complex of dimension at most $d-1$ vanishes
above degree $d-1$. Hence $\widetilde H_i(\X;\mathbb Z)=0$ for
$i\ge d$, and $\widetilde H_{d-1}(\X;\mathbb Z)$ is the kernel of the
boundary map on the group of $(d-1)$-chains, which is a subgroup of a
finitely generated free abelian group and hence free.
\end{proof}

\subsection{Stanley--Reisner rings}
\label{subsec:SR}

Let $S=\Bbbk[x_v:v\in V]$. The \emph{Stanley--Reisner ideal} $I_\X$
of $\X$ is the ideal of $S$ generated by the monomials
$x_\tau=\prod_{v\in\tau}x_v$, where $\tau$ ranges over the minimal
non-faces of $\X$. By Hochster's formula (see, e.g.,~\cite{MS}),
\begin{equation}\label{eq:hochster-betti}
  \beta_{i,W}(I_\X)=\dim_\Bbbk\widetilde H_{|W|-i-2}(\X[W];\Bbbk)
  ~(W\subseteq V,\ i\ge 0).
\end{equation}
Consequently, $\reg(S/I_\X)=\ly(\X)$, and $\reg I_\X=\ly(\X)+1$
whenever $I_\X\ne 0$. For a graph $G$, the Stanley--Reisner ideal of
$\Ind(G)$ is the \emph{edge ideal} $I(G)=(x_ux_v:uv\in E(G))$. Thus
$\reg I(G)=\ly(\Ind(G))+1$ whenever $E(G)\ne\emptyset$.

\subsection{Co-chordal covers}
\label{subsec:cochord}

A graph is \emph{chordal} if it has no induced cycle of length at
least four, and \emph{co-chordal} if its complement is chordal. A
vertex is \emph{simplicial} if its neighborhood is a clique. By
results of Dirac~\cite{Dir} and Fulkerson and Gross~\cite{FG}, a
graph is chordal if and only if it admits a \emph{perfect elimination
ordering}, that is, an ordering of its vertices in which each vertex
is simplicial in the subgraph induced by itself and its successors.
The \emph{co-chordal cover number} $\cochord(G)$ is the minimum number
of co-chordal subgraphs of $G$ whose edge sets cover
$E(G)$~\cite{W14}. These subgraphs may be assumed to be spanning,
since adding isolated vertices to a graph adds dominating vertices to
its complement and therefore preserves co-chordality.

Woodroofe~\cite{W14} proved that $\ly(\Ind(G))\le\cochord(G)$, or
equivalently $\reg I(G)\le\cochord(G)+1$. The following theorem
strengthens this bound to the $\theta$-number, and hence to the
collapsibility number.

\begin{theorem}\label{thm:ind-cochord}
For every graph $G$,
\begin{equation}\label{eq:ind-cochord}
  \im(G)\le\ly(\Ind(G))\le\coll(\Ind(G))\le\theta(\Ind(G))
  \le\cochord(G).
\end{equation}
\end{theorem}

\begin{proof}
If $G$ has no edges, then $\Ind(G)$ is a simplex and every term
in~\eqref{eq:ind-cochord} equals $0$. Assume therefore that $G$ has
an edge.
Let $\{e_1,\ldots,e_m\}$ be an induced matching of $G$, and let
$U=\bigcup_{i=1}^m e_i$. Then
\[
  \Ind(G)[U]=\Ind(G[e_1])*\cdots*\Ind(G[e_m])
\]
is the boundary of an $m$-dimensional cross-polytope, which is a
sphere of dimension $m-1$. Its reduced homology in degree $m-1$ is
nonzero over every field, and hence $\ly(\Ind(G))\ge m$; see
also~\cite{Kat}. The two middle inequalities are~\ref{F:leray}.
For the last inequality, recall that $\theta(\Ind(H))\le 1$ for every
co-chordal graph $H$~\cite[Proposition~19]{BC}. Let $c=\cochord(G)$,
and let $H_1,\ldots,H_c$ be spanning co-chordal subgraphs of $G$
whose edge sets cover $E(G)$. A subset of $V(G)$ is independent in
$G$ if and only if it is independent in every $H_j$, so
$\Ind(G)=\bigcap_{j=1}^c\Ind(H_j)$. By
Proposition~\ref{prop:theta-subadditivity},
$\theta(\Ind(G))\le\sum_{j=1}^c\theta(\Ind(H_j))\le c.$
This completes the proof.
\end{proof}

\begin{remark}
If $\im(G)=\cochord(G)$, then Theorem~\ref{thm:ind-cochord} yields
\[
  \im(G)=\ly(\Ind(G))=\coll(\Ind(G))=\theta(\Ind(G))=\cochord(G)
\]
over every field. This is the case for weakly chordal graphs, by a
result of Busch, Dragan, and Sritharan~\cite[Proposition~3]{BDS},
and for well-covered bipartite graphs, by a result of
Woodroofe~\cite[Theorem~16]{W14}. It is also the case whenever
$\im(G)=\nu(G)$, where $\nu(G)$ denotes the matching number of $G$,
because $\cochord(G)\le\nu(G)$. Indeed, if $M$ is a maximum
matching, then for each edge $ab\in M$ the edges of $G$ meeting
$\{a,b\}$ form a co-chordal subgraph, and these $|M|$ subgraphs
cover $E(G)$ by the maximality of $M$.
\end{remark}

\subsection{Lew's invariant}
\label{subsec:lew}

Lew~\cite{Lew} introduced the following non-recursive upper bound for
the collapsibility number.

\begin{theorem}[{\cite[Theorem~4]{Lew}}]\label{thm:lew}
Let $\mathcal S(\X)$ be the family of sets $T$ of vertices of $\X$
for which there exist an enumeration $v_1,\ldots,v_k$ of $T$ and
facets $\sigma_1,\ldots,\sigma_{k+1}$ of $\X$ such that
$v_i\notin\sigma_i$ for $1\le i\le k$ and $v_i\in\sigma_j$ for
$1\le i<j\le k+1$. Let $d'(\X)=\max\{|T|:T\in\mathcal S(\X)\}$,
where $T=\emptyset$ is allowed. Then $\coll(\X)\le d'(\X)$.
\end{theorem}

For independence complexes, Lew obtained the following
graph-theoretic description of $d'$.

\begin{lemma}[{\cite[Proposition~11]{Lew}}]\label{lem:dprime-flag}
Let $G$ be a graph. Then $d'(\Ind(G))$ is the largest integer $k$
for which there exist vertices
$p_1,\ldots,p_k,q_1,\ldots,q_k\in V(G)$ such that
\begin{enumerate}[label=\textup{(\roman*)}]
\item $p_1,\ldots,p_k$ are distinct and pairwise nonadjacent;
\item $p_iq_i\in E(G)$ for $1\le i\le k$;
\item $p_jq_i\notin E(G)$ for $1\le j<i\le k$.
\end{enumerate}
\end{lemma}

We call such a sequence $(p_i,q_i)_{i=1}^k$ a \emph{Lew sequence},
and we call $q_i$ the \emph{witness} of $p_i$. The $2k$ vertices of a
Lew sequence are distinct. Indeed, $q_i\ne p_i$ by~(ii); if
$q_i=p_j$ with $j\ne i$, then $p_j$ and $p_i$ are adjacent by~(ii),
contrary to~(i); and if $q_i=q_j$ with $j<i$, then (ii) and (iii)
contradict each other.

\begin{lemma}\label{lem:dprime-dim}
For every simplicial complex $\X$, we have $d'(\X)\le\dim\X+1$.
\end{lemma}

\begin{proof}
In the notation of Theorem~\ref{thm:lew}, every $T\in\mathcal S(\X)$
is contained in the facet $\sigma_{k+1}$. Hence
$|T|\le|\sigma_{k+1}|\le\dim\X+1$.
\end{proof}

\begin{theorem}\label{thm:theta-dprime}
For every simplicial complex $\X$,
\[
  \coll(\X)\le\theta(\X)\le d'(\X)\le\dim\X+1.
\]
\end{theorem}

\begin{proof}
The first inequality is part of~\ref{F:leray}, and the third is
Lemma~\ref{lem:dprime-dim}. For the second, B{\i}y{\i}ko\u{g}lu and
Civan~\cite[Corollary~29]{BC} prove that $\theta(\X)\le k(\X)$, where
$k(\X)$ is the maximum size of a set in a family that coincides with
the family $\mathcal S(\X)$ of Theorem~\ref{thm:lew}. Hence
$k(\X)=d'(\X)$, and the second inequality follows.
\end{proof}

\begin{remark}
The inequality $\theta(\X)\le d'(\X)$ also follows by combining the
bound $\theta(\X)\le d(\X,\prec)$ of Santhanam, Shukla, and
Singh~\cite[Proposition 3.19]{SSS}, valid for every linear order $\prec$ on the facets
of $\X$, with Lew's inequality $d(\X,\prec)\le d'(\X)$~\cite{Lew}.
\end{remark}

\begin{proof}[Proof of Theorem~\textup{\ref{thmB}}]
The chain of inequalities follows from
Theorems~\ref{thm:ind-cochord} and~\ref{thm:theta-dprime}, and the
bound $d'(\Ind(G))\le\dim\Ind(G)+1$ is Lemma~\ref{lem:dprime-dim}.
For the homotopy statement, suppose that $G$ has an edge, and let
$d=\theta(\Ind(G))$. Then $d\ge\im(G)\ge 1$, so
Lemma~\ref{lem:theta-homotopy} shows that $\Ind(G)$ is homotopy
equivalent to a finite CW complex of dimension at most
$d-1\le\min\{d'(\Ind(G)),\cochord(G)\}-1$.
\end{proof}

\section{Vietoris--Rips complexes of hypercube graphs}
\label{sec:VR}

\subsection{Setup}
\label{subsec:VR-setup}

We use the notation of Section~\ref{sec:intro}. For
$\alpha\in\F^n$, we write $\wt(\alpha)=|\supp(\alpha)|$ for the
Hamming weight of $\alpha$, so that
$d_H(\alpha,\beta)=\wt(\alpha+\beta)$. For $Y\subseteq\F^n$, we have
$\VR(Y;r)=\VR(I_n;r)[Y]$. Let $e_1,\ldots,e_n$ denote the standard
unit vectors, let $\one=(1,\ldots,1)$, and for $\alpha\in\F^n$ let
$\overline{\alpha}=\alpha+\one$ denote the antipode of $\alpha$.
Binary strings are concatenated by juxtaposition. A
\emph{$k$-dimensional coordinate subcube} of $\F^n$ is a set obtained
by fixing $n-k$ coordinates and letting the remaining $k$ coordinates
vary.

Recall that two vertices of $H_{n,r}$ are adjacent if and only if
their Hamming distance exceeds $r$; equivalently,
$H_{n,r}=\Cay(\F^n,\{\gamma:\wt(\gamma)>r\})$. Translations and
coordinate permutations act on $H_{n,r}$ by automorphisms, and the
translations alone act transitively on the vertices, so $H_{n,r}$ is
vertex-transitive. A $k$-dimensional coordinate subcube is isometric
to $\F^k$, since its elements agree in the fixed coordinates. Hence
the subgraph of $H_{n,r}$ that it induces is isomorphic to $H_{k,r}$,
which is edgeless when $r\ge k$.

Throughout this section, $\Bbbk$ is a field and
$S=\Bbbk[x_\alpha:\alpha\in\F^n]$. Unless stated otherwise, Leray
numbers and regularities are taken over $\Bbbk$.

\begin{proposition}\label{prop:VR-flag}
Let $0\le r\le n$, and let $\Delta=\VR(I_n;r)$. Then $\Delta$ is a
flag complex, $\Delta=\Ind(H_{n,r})$, and
$I_\Delta=I(H_{n,r})=(x_\alpha x_\beta:d_H(\alpha,\beta)>r).$
If $r<n$, then
\[
  \reg(S/I_\Delta)=\ly(\Delta)\le\coll(\Delta)\le\theta(\Delta)
  \le\min\{d'(\Delta),\cochord(H_{n,r})\}
\]
and $\reg I_\Delta=\ly(\Delta)+1$.
\end{proposition}

\begin{proof}
A subset of $\F^n$ has diameter at most $r$ if and only if no two of
its elements are at Hamming distance greater than $r$, that is, if
and only if it is independent in $H_{n,r}$. Thus
$\Delta=\Ind(H_{n,r})$; in particular, $\Delta$ is flag, and the
description of $I_\Delta$ follows. The inequalities follow
from~\eqref{eq:hochster-betti} and
Theorem~\ref{thmB}. If $r<n$, then antipodal vertices
are at distance $n>r$, so $H_{n,r}$ has an edge. Hence $I_\Delta\ne 0$
and $\reg I_\Delta=\ly(\Delta)+1$.
\end{proof}

If $r\ge n$, then $\VR(I_n;r)$ is a simplex. If $r=n-1$, then two
vertices are adjacent in $H_{n,n-1}$ if and only if they are at
distance $n$, that is, antipodal. Thus $H_{n,n-1}$ is the perfect
matching joining antipodal vertices, and $\VR(I_n;n-1)$ is the
boundary of a $2^{n-1}$-dimensional cross-polytope~\cite{AA}. A graph
$G$ that is a matching with $m$ edges satisfies $\im(G)=m$, and also
$\cochord(G)\le m$, since each edge, together with the remaining
vertices as isolated vertices, forms a co-chordal subgraph. Hence
$\im(G)=\cochord(G)=m$, and Theorem~\ref{thm:ind-cochord} yields
\begin{equation}\label{eq:VR-top}
  \ly(\VR(I_n;n-1))=\coll(\VR(I_n;n-1))=\theta(\VR(I_n;n-1))
  =2^{n-1}.
\end{equation}

\subsection{The induced matching number}
\label{subsec:VR-matching}

We first obtain a lower bound from coordinate subcubes.

\begin{theorem}\label{thm:VR-matching}
Let $0\le r<n$, and let $Q\subseteq\F^n$ be an $(r+1)$-dimensional
coordinate subcube. Then $H_{n,r}[Q]$ is a perfect matching with
$2^r$ edges, and $\VR(Q;r)$ is the boundary of a $2^r$-dimensional
cross-polytope; in particular, $\VR(Q;r)$ is homeomorphic to the
sphere $S^{2^r-1}$. Consequently, $\im(H_{n,r})\ge 2^r$, and for
every field $\Bbbk$,
\[
  2^r\le\reg(S/I_{\VR(I_n;r)})=\ly(\VR(I_n;r))
  \le\coll(\VR(I_n;r))\le\theta(\VR(I_n;r)).
\]
\end{theorem}

\begin{proof}
The subcube $Q$ is isometric to $\F^{r+1}$, in which two vertices are
at Hamming distance greater than $r$ if and only if they are
antipodal. Thus $H_{n,r}[Q]$ is a perfect matching on $2^{r+1}$
vertices, with $2^r$ edges. Consequently,
$\VR(Q;r)=\Ind(H_{n,r}[Q])$ is the join of $2^r$ copies of $S^0$,
which is the boundary of a $2^r$-dimensional cross-polytope. The
remaining assertions follow from Theorem~\ref{thm:ind-cochord} and
Proposition~\ref{prop:VR-flag}.
\end{proof}

The matching in Theorem~\ref{thm:VR-matching} is in fact a maximum
induced matching. To prove this, we use the following two-sided
Hamming set-pair inequality; see also~\cite[Section~4.3]{AJS} for a
related application to induced matchings.

\begin{theorem}[Alon, Jin, and Sudakov~{\cite[Theorem~1.5]{AJS}}]
\label{thm:AJS}
Let $A$ be an alphabet with $|A|\ge 2$, and let $n>t\ge 0$. Suppose
that $a_1,\ldots,a_m,b_1,\ldots,b_m\in A^n$ satisfy
$d_H(a_i,b_i)\ge t+1$ for every $i$ and $d_H(a_i,b_j)\le t$ whenever
$i\ne j$. Then $m\le 2^{t+1}$.
\end{theorem}

\begin{corollary}\label{thm:im-exact}
For all integers $n>r\ge 0$, we have $\im(H_{n,r})=2^r$.
\end{corollary}

\begin{proof}
The lower bound $\im(H_{n,r})\ge 2^r$ is part of
Theorem~\ref{thm:VR-matching}. For the upper bound, let
$u_1y_1,\ldots,u_my_m$ be an induced matching in $H_{n,r}$, and
consider the $2m$ ordered pairs $(u_i,y_i)$ and $(y_i,u_i)$ for
$1\le i\le m$. The two entries of each pair are at Hamming distance
greater than $r$, since $u_iy_i$ is an edge of $H_{n,r}$. Let $(a,b)$
and $(a',b')$ be two distinct pairs. If they arise from different
edges of the matching, then $d_H(a,b')\le r$, because the matching is
induced. If they arise from the same edge, then $a=b'$, so
$d_H(a,b')=0$. Thus the $2m$ pairs satisfy the hypotheses of
Theorem~\ref{thm:AJS} with $A=\F$ and $t=r$, and therefore
$2m\le 2^{r+1}$, that is, $m\le 2^r$.
\end{proof}

As an immediate consequence of Corollary~\ref{thm:im-exact}, we record
the following colorful diameter statement for the hypercube.

\begin{corollary}[colorful diameter theorem]
\label{cor:colourful-diameter}
Let $n>r\ge 1$, and let $V_1,\ldots,V_{2^r+1}$ be nonempty subsets of
$\F^n$. Suppose that for every choice of points $v_i\in V_i$,
$1\le i\le 2^r+1$, the set $\{v_1,\ldots,v_{2^r+1}\}$ has Hamming
diameter at most $r$. Then some $V_i$ has Hamming diameter at most
$r$. The number $2^r+1$ cannot be replaced by $2^r$, even if the sets
$V_i$ are required to be pairwise disjoint.
\end{corollary}

\begin{proof}
Suppose, to the contrary, that every $V_i$ contains points $a_i$ and
$b_i$ with $d_H(a_i,b_i)>r$. Let $i\ne j$, let
$x\in\{a_i,b_i\}$, and let $y\in\{a_j,b_j\}$. Since every $V_k$ is
nonempty, there is a choice of points $v_k\in V_k$ with $v_i=x$ and
$v_j=y$, and the hypothesis gives $d_H(x,y)\le r$. In particular, the
pairs $\{a_i,b_i\}$ are pairwise disjoint: if $x\in\{a_i,b_i\}$ also
belonged to $\{a_j,b_j\}$, then applying the previous observation to
$x$ and the other point $x'$ of $\{a_i,b_i\}$ would give
$d_H(x,x')\le r$, which is a contradiction. Hence the edges
$a_ib_i$, $1\le i\le 2^r+1$, form an induced matching of size
$2^r+1$ in $H_{n,r}$, contradicting Corollary~\ref{thm:im-exact}.
To prove sharpness, let $Q$ be an $(r+1)$-dimensional coordinate
subcube, and let $V_1,\ldots,V_{2^r}$ be its $2^r$ antipodal pairs.
Each $V_i$ has diameter $r+1$. On the other hand, two points chosen
from different pairs are not antipodal in $Q$, so they are at
distance at most $r$. Thus every choice of one point from each $V_i$
has diameter at most $r$.
\end{proof}

\subsection{A strengthening of Shukla's conjecture}
\label{subsec:VR-conj}

The exact values in Theorem~\ref{thmA} suggest the following
strengthening of Conjecture~\ref{conj:shukla}.

\begin{conjecture}\label{conj:theta}
For all integers $n$ and $r$ with $n\ge r+1\ge 2$,
$\theta(\VR(I_n;r))=2^r.$
\end{conjecture}

Since $\coll\le\theta$ by~\ref{F:leray}, and $\coll(\VR(I_n;r))\ge 2^r$
by Theorem~\ref{thm:VR-matching}, Conjecture~\ref{conj:theta} implies
Conjecture~\ref{conj:shukla}.

Conjecture~\ref{conj:theta} can be restated in terms of arbitrary
subsets of the hypercube. Indeed, it is equivalent to the assertion
that $\theta(\VR(Y;r))\le 2^r$ for all integers $r\ge 1$ and
$n\ge 1$ and every nonempty $Y\subseteq\F^n$. If $n\ge r+1$, this
follows from Conjecture~\ref{conj:theta} by~\ref{F:mono}; if
$n\le r$, then $\VR(Y;r)$ is a simplex and $\theta(\VR(Y;r))=0$.
Conversely, taking $Y=\F^n$ and using the lower bound in
Theorem~\ref{thm:VR-matching} recovers Conjecture~\ref{conj:theta}.
By Lemma~\ref{lem:theta-homotopy}, Conjecture~\ref{conj:theta} would
imply that every such complex $\VR(Y;r)$ is homotopy equivalent to a
finite CW complex of dimension at most $2^r-1$, and hence that
$\widetilde H_i(\VR(Y;r);\mathbb Z)=0$ for all $i\ge 2^r$. In
particular, it would imply the vanishing of homology predicted by
Adamaszek and Adams~\cite[Question~9]{AA}, although not necessarily
the collapse onto the $(2^r-1)$-skeleton asked for there. Moreover,
combined with Theorem~\ref{thm:VR-matching} and Hochster's formula,
it would give $\reg I_{\VR(I_n;r)}=2^r+1$ for all $n>r\ge 1$ over
every field.

\begin{remark}\label{rem:betti-short}
By Hochster's formula~\eqref{eq:hochster-betti}, every
$(r+1)$-dimensional coordinate subcube $Q$ contributes
$\dim_\Bbbk\widetilde H_{2^r-1}(\VR(Q;r);\Bbbk)=1$ to the graded Betti
number $\beta_{2^r-1,\,2^{r+1}}(I_{\VR(I_n;r)})$. Hence this Betti
number is at least $2^{n-r-1}\binom{n}{r+1}$ over every field, which
is the Betti-number analogue of the homological lower bound of Adams
and Virk~\cite[Theorem~4.1]{AV}; unlike the homological bound, it
also holds for $r=1$.
\end{remark}
\subsection{Monotonicity in the dimension}
\label{subsec:VR-mono}

We next show that the invariants are monotone in the dimension of the
ambient hypercube.

\begin{proposition}\label{prop:VR-mono}
Fix $r\ge 1$. For $n>r$, let $f(n)$ denote any one of
\[
  \im(H_{n,r}),~\ly(\VR(I_n;r)),~\coll(\VR(I_n;r)),~
  \theta(\VR(I_n;r)),~ d'(\VR(I_n;r)),~\cochord(H_{n,r}).
\]
Then $f(n)\le f(n+1)$. In particular, if $\theta(\VR(I_N;r))\le 2^r$
for some $N>r$, then
\[
  \ly(\VR(I_n;r))=\coll(\VR(I_n;r))=\theta(\VR(I_n;r))=2^r
  ~(r<n\le N).
\]
\end{proposition}

\begin{proof}
The coordinate slice $\{\alpha0:\alpha\in\F^n\}\subseteq\F^{n+1}$ is
isometric to $\F^n$. Hence $H_{n,r}$ is an induced subgraph of
$H_{n+1,r}$, and $\VR(I_n;r)$ is an induced subcomplex of
$\VR(I_{n+1};r)$. The monotonicity of $\ly$, $\coll$, and $\theta$
follows from~\ref{F:mono}. A Lew sequence in an induced subgraph is
a Lew sequence in the whole graph, so the monotonicity of $d'$
follows from Lemma~\ref{lem:dprime-flag}. Similarly, an induced
matching of an induced subgraph is an induced matching of the whole
graph. Finally, if $H_1,\ldots,H_c$ are co-chordal subgraphs of
$H_{n+1,r}$ whose edge sets cover its edges, then their restrictions
to the slice are co-chordal and cover the edges of $H_{n,r}$. This
proves that $f(n)\le f(n+1)$ in all six cases. The final assertion
follows from this monotonicity, Theorem~\ref{thm:VR-matching},
and~\ref{F:leray}.
\end{proof}

\subsection{Vertex links}
\label{subsec:VR-links}

For $v\in\F^n$, let $B(v,r)=\{\alpha\in\F^n:d_H(\alpha,v)\le r\}$,
and let $B^*_{n,r}=B(0,r)\setminus\{0\}$. Since $\VR(I_n;r)$ is a
flag complex, the link of a vertex $v$ is the subcomplex induced on
its neighbors in the $1$-skeleton. Combined with translation by $v$,
this gives
\[
  \lk(\VR(I_n;r),v)=\VR(B(v,r)\setminus\{v\};r)
  \cong\VR(B^*_{n,r};r).
\]

\begin{proposition}\label{prop:VR-link}
Let $1\le r<n$, and let $\Delta=\VR(I_n;r)$ and
$P=\VR(B^*_{n,r};r)$. Then
\[
\begin{aligned}
  2^r-1&\le\ly(P)\le\ly(\Delta)\le\ly(P)+1,\\
  2^r-1&\le\theta(P)\le\theta(\Delta)\le\theta(P)+1.
\end{aligned}
\]
Consequently, the following hold.
\begin{enumerate}[label=\textup{(\alph*)}]
\item If $\theta(P)=2^r-1$, then
$\ly(\Delta)=\coll(\Delta)=\theta(\Delta)=2^r$. In particular,
Conjecture~\ref{conj:theta} holds for $(n,r)$.
\item If $\theta(\Delta)=2^r$, then $\theta(P)\in\{2^r-1,2^r\}$.
\end{enumerate}
\end{proposition}

\begin{proof}
Choose an $(r+1)$-dimensional coordinate subcube $Q$ containing $0$,
and let $u$ be the antipode of $0$ in $Q$. Every element of $Q$
other than $u$ lies within distance $r$ of $0$, so
$W=Q\cap B^*_{n,r}=Q\setminus\{0,u\}$. By
Theorem~\ref{thm:VR-matching}, $H_{n,r}[Q]$ is a perfect matching
with $2^r$ edges, one of which is $0u$. Hence $H_{n,r}[W]$ is a
perfect matching with $2^r-1$ edges, and $\VR(W;r)$ is the boundary
of a $(2^r-1)$-dimensional cross-polytope. Since $\VR(W;r)$ is an
induced subcomplex of $P$, it follows that $\ly(P)\ge 2^r-1$, and
hence $\theta(P)\ge 2^r-1$ by~\ref{F:leray}.
Since $P$ is an induced subcomplex of $\Delta$,~\ref{F:mono} gives
$\ly(P)\le\ly(\Delta)$ and $\theta(P)\le\theta(\Delta)$. Moreover,
$\Delta$ is not a simplex because $r<n$, and translations act
transitively on its vertices, so every vertex link is isomorphic to
$P$. Lemma~\ref{lem:link} therefore yields
$\ly(\Delta)\le\ly(P)+1$ and $\theta(\Delta)\le\theta(P)+1$. This
proves both chains of inequalities.
If $\theta(P)=2^r-1$, then $\theta(\Delta)\le 2^r$, while
Theorem~\ref{thm:VR-matching} and~\ref{F:leray} give
$2^r\le\ly(\Delta)\le\coll(\Delta)\le\theta(\Delta)$; this
proves~(a). Part~(b) follows directly from the second chain.
\end{proof}

Proposition~\ref{prop:VR-link} shows that the condition
$\theta(\VR(B^*_{n,r};r))=2^r-1$ is sufficient for
Conjecture~\ref{conj:theta} to hold for $(n,r)$. For fixed $r$, this
replaces the $2^n$ vertices of the hypercube by the
$\sum_{i=1}^r\binom{n}{i}$ vertices of the punctured Hamming ball. We
do not know whether the condition is also necessary. We verify it in
the parameter ranges of Theorem~\ref{thmA}
(Lemma~\ref{lem:punctured}), and we use the corresponding bounds for
links again in Section~\ref{sec:topology}.

\section{Exact values: small scales and small codimension}
\label{sec:small}

In this section we prove Theorem~\ref{thmA}. We treat the scale
$r=1$ by an explicit cover by two co-chordal subgraphs, the scale
$r=2$ by a bound on Lew's invariant in terms of common
non-neighborhoods of edges, and the case $n=r+2$ by an explicit
co-chordal cover.

\subsection{Small scales}
\label{subsec:small-scales}

\begin{lemma}\label{lem:common}
Let $G$ be a graph with at least one edge, and let
\[
  \mu(G)=\max_{ab\in E(G)}\alpha\bigl(G-(N[a]\cup N[b])\bigr).
\]
Then $\theta(\Ind(G))\le d'(\Ind(G))\le\mu(G)+1$. If $G$ is a perfect matching, that is, a disjoint union of $m\ge 1$
edges, then both inequalities are equalities; more precisely,
$\theta(\Ind(G))=d'(\Ind(G))=\mu(G)+1=m$.
\end{lemma}

\begin{proof}
The first inequality is part of Theorem~\ref{thm:theta-dprime}. For
the second, let $k=d'(\Ind(G))$, which is at least $1$ because $G$
has an edge, and let $(p_i,q_i)_{i=1}^k$ be a Lew sequence of length
$k$. For $j<k$, the vertex $p_j$ is not in $N[p_k]$ by
condition~(i) of Lemma~\ref{lem:dprime-flag}, and it is not in
$N[q_k]$ by condition~(iii) and because the vertices of a Lew
sequence are distinct. Hence $\{p_1,\ldots,p_{k-1}\}$ is an
independent set in $G-(N[p_k]\cup N[q_k])$. Since $p_kq_k\in E(G)$,
it follows that $k-1\le\mu(G)$.
Now suppose that $G$ is a perfect matching with $m$ edges. For every
edge $ab$, we have $N[a]\cup N[b]=\{a,b\}$, so
$G-(N[a]\cup N[b])$ is a perfect matching with $m-1$ edges, and
$\mu(G)=m-1$. On the other hand, $\theta(\Ind(G))=m$ by
Theorem~\ref{thm:ind-cochord}, since $\im(G)=\cochord(G)=m$. Hence
$m\le\theta(\Ind(G))\le d'(\Ind(G))\le m$.
\end{proof}

\begin{proposition}\label{prop:VR-r1}
For every $n\ge 2$, we have $\im(H_{n,1})=\cochord(H_{n,1})=2$ and
\[
  \ly(\VR(I_n;1))=\coll(\VR(I_n;1))=\theta(\VR(I_n;1))
  =d'(\VR(I_n;1))=2.
\]
\end{proposition}

\begin{proof}
Since $I_n$ is bipartite, it has no triangles, so $\VR(I_n;1)=I_n$ is
one-dimensional. Hence $d'(\VR(I_n;1))\le 2$ by
Lemma~\ref{lem:dprime-dim}, while $\im(H_{n,1})\ge 2$ by
Theorem~\ref{thm:VR-matching}.

Let $A_0$ and $A_1$ denote the sets of vertices of $\F^n$ of even and
odd weight, respectively. For $j\in\{0,1\}$, let $P_j$ be the
spanning subgraph of $H_{n,1}$ consisting of all edges that meet
$A_j$. No edge of $P_j$ has both endpoints in $A_{1-j}$, so $A_{1-j}$
is a clique in $\overline{P_j}$. Two distinct vertices of $A_j$ are at
even Hamming distance, hence at distance at least $2$, so they are adjacent
in $H_{n,1}$ and therefore in $P_j$; thus $A_j$ is independent in
$\overline{P_j}$. Consequently, $\overline{P_j}$ is a split graph and
is therefore chordal. Since every edge of $H_{n,1}$ meets $A_0$ or
$A_1$, the graphs $P_0$ and $P_1$ cover $E(H_{n,1})$, so
$\cochord(H_{n,1})\le 2$. All the claimed equalities now follow from
Theorem~\ref{thm:ind-cochord} and Theorem~\ref{thmB}.
\end{proof}

At scale $2$, the bound of Lemma~\ref{lem:common} determines Lew's
invariant.

\begin{theorem}\label{thm:VR-r2}
For every $n\ge 1$ and every $Y\subseteq\F^n$, we have
$\theta(\VR(Y;2))\le 4$. Moreover, if $n\ge 3$, then over every field
\[
\begin{aligned}
  4&=\im(H_{n,2})=\reg(S/I_{\VR(I_n;2)})=\ly(\VR(I_n;2))\\
   &=\coll(\VR(I_n;2))=\theta(\VR(I_n;2))=d'(\VR(I_n;2)).
\end{aligned}
\]
Finally, every complex $\VR(Y;2)$ is homotopy equivalent to a finite
CW complex of dimension at most $3$.
\end{theorem}

\begin{proof}
By~\ref{F:mono}, it suffices to prove that
$\theta(\VR(I_n;2))\le 4$. If $n\le 2$, then $\VR(I_n;2)$ is a
simplex, so we may assume that $n\ge 3$. Let $H=H_{n,2}$. We show that
$\mu(H)=3$.
Fix an edge $ab$ of $H$. By translation, we may assume that $a=0$.
Let $D=\supp(b)$ and $s=|D|$, so that $s\ge 3$. Let
$W=V(H)\setminus(N[a]\cup N[b])$. A vertex $u$ belongs to $W$ if and
only if $d_H(u,0)\le 2$ and $d_H(u,b)\le 2$; these conditions exclude
$u=0$ and $u=b$, since $d_H(0,b)=s>2$. Write
$\supp(u)=X_1\sqcup X_2$ with $X_1\subseteq D$ and
$X_2\cap D=\emptyset$. The two conditions become
\[
  |X_1|+|X_2|\le 2~
  \text{and}~
  s-|X_1|+|X_2|\le 2,
\]
and adding them gives $s+2|X_2|\le 4$. If $s\ge 5$, then
$W=\emptyset$. If $s\in\{3,4\}$, then $X_2=\emptyset$ and
$s-2\le|X_1|\le 2$. Hence $W$ consists of the six nonempty proper
subsets of $D$ if $s=3$, and of the six $2$-element subsets of $D$
if $s=4$. In either case, two vertices of $W$ are at Hamming distance
greater than $2$ if and only if their supports are complementary in
$D$. Thus $H[W]$ is a matching with three edges, and
$\alpha(H[W])=3$. Since $n\ge 3$, edges with $s=3$ exist, and hence
$\mu(H)=3$.
Lemma~\ref{lem:common} now gives
$\theta(\VR(I_n;2))\le d'(\VR(I_n;2))\le 4$. The reverse inequalities
follow from Theorems~\ref{thm:VR-matching} and~\ref{thm:ind-cochord},
together with Proposition~\ref{prop:VR-flag}. Finally, the homotopy
statement follows from the bound $\theta(\VR(Y;2))\le 4$ and
Lemma~\ref{lem:theta-homotopy}.
\end{proof}

\begin{remark}\label{rem:r2}
Theorem~\ref{thm:VR-r2} gives a new proof of Shukla's equality
$\coll(\VR(I_n;2))=4$~\cite{Shu}, and it also determines the
$\theta$-number and Lew's invariant. The uniform bound on the
homotopy dimension is consistent with the theorem of Adamaszek and
Adams~\cite{AA} that $\VR(I_n;2)$ is homotopy equivalent to a wedge
of $3$-spheres. By Kleitman's diametric theorem~\cite{Kle},
$\dim\VR(I_n;2)=n$ for $n\ge 3$, the maximum being attained by
Hamming balls of radius one. Thus the homotopy dimension remains
bounded, although the dimension of the complex grows with $n$.
\end{remark}

\subsection{Codimension two}
\label{subsec:codim-two}

We prove Shukla's conjecture for $n=r+2$ by constructing an optimal
co-chordal cover. We first give a general criterion for a subgraph to
be co-chordal.

\begin{lemma}\label{lem:star-clique}
Let $G$ be a graph, let $K$ be a clique of $G$, and let
$F\subseteq E(G-K)$. Suppose that the graph $(V(G)\setminus K,F)$ is
co-chordal and that, for every $c\in K$ and every $ab\in F$, at least
one of $a$ and $b$ is adjacent to $c$ in $G$. Then the spanning
subgraph $P$ of $G$ with edge set
\[
  E(P)=F\cup\{ab\in E(G):\{a,b\}\cap K\ne\emptyset\}
\]
is co-chordal.
\end{lemma}

\begin{proof}
We show that $\overline P$ is chordal. Let $c\in K$. Since no edge
of $F$ meets $K$, we have $N_P(c)=N_G(c)$, and hence
$N_{\overline P}(c)=V(G)\setminus N_G[c]$. Since $K$ is a clique,
this set is disjoint from $K$. The edges of $P$ joining two vertices
outside $K$ are exactly the edges of $F$, and by hypothesis no edge
of $F$ joins two non-neighbors of $c$. Hence $N_{\overline P}(c)$ is
a clique of $\overline P$, that is, $c$ is simplicial in
$\overline P$. A simplicial vertex remains simplicial in every
induced subgraph containing it, so the vertices of $K$ can be
eliminated first, in any order. The remaining graph is the
complement of $(V(G)\setminus K,F)$, which is chordal by hypothesis.
Hence $\overline P$ admits a perfect elimination ordering.
\end{proof}

\begin{theorem}\label{thm:VR-n-2}
For every $r\ge 1$, we have $\im(H_{r+2,r})=\cochord(H_{r+2,r})=2^r$.
Consequently, over every field,
\[
\begin{aligned}
  \reg(S/I_{\VR(I_{r+2};r)})&=\ly(\VR(I_{r+2};r))
  =\coll(\VR(I_{r+2};r))
  =\theta(\VR(I_{r+2};r))=2^r.
\end{aligned}
\]
Moreover, $\VR(I_{r+2};r)$ is homotopy equivalent to a finite CW
complex of dimension at most $2^r-1$. In particular,
Conjecture~\ref{conj:shukla} holds for $n=r+2$.
\end{theorem}

\begin{proof}
Let $H=H_{r+2,r}$, and write each vertex of $H$ as $yt$ with
$y\in\F^{r+1}$ and $t\in\F$. Let $\mathcal P$ denote the set of the
$2^r$ antipodal pairs in $\F^{r+1}$. For $p=\{y,\overline y\}\in
\mathcal P$, let
\[
  K_p=\{y0,\overline y0\}
  ~\text{and}~
  F_p=\{\{y1,\overline y1\}\},
\]
and let $P_p$ be the spanning subgraph of $H$ whose edges are the
edge in $F_p$ together with all edges of $H$ that meet $K_p$. Since
$d_H(y0,\overline y0)=r+1$, the set $K_p$ is a clique of $H$; since
$d_H(y1,\overline y1)=r+1$, the edge in $F_p$ is an edge of $H-K_p$.
The graph $(V(H)\setminus K_p,F_p)$ has a single edge and is
therefore co-chordal. Moreover,
$d_H(y0,\overline y1)=d_H(\overline y0,y1)=r+2$, so each vertex of
$K_p$ is adjacent to an endpoint of the edge in $F_p$. By
Lemma~\ref{lem:star-clique}, $P_p$ is co-chordal.

We claim that the graphs $P_p$, $p\in\mathcal P$, cover $E(H)$. The
sets $K_p$ partition the layer $\F^{r+1}\times\{0\}$, so every edge
of $H$ meeting this layer belongs to some $P_p$. Two vertices $y1$
and $z1$ of the layer $\F^{r+1}\times\{1\}$ are adjacent in $H$ if
and only if $d_H(y,z)=r+1$, that is, if and only if $y$ and $z$ are
antipodal in $\F^{r+1}$; these edges are exactly the edges in the
sets $F_p$. This proves the claim, and hence $\cochord(H)\le 2^r$.
Since $\im(H)=2^r$ by Corollary~\ref{thm:im-exact},
Theorem~\ref{thm:ind-cochord} gives $\cochord(H)\ge 2^r$, and hence
$\im(H)=\cochord(H)=2^r$. The remaining equalities follow from
Theorem~\ref{thm:ind-cochord} and Proposition~\ref{prop:VR-flag},
and the homotopy statement follows from
Lemma~\ref{lem:theta-homotopy}.
\end{proof}

\subsection{Punctured Hamming balls}
\label{subsec:punctured}

The following local counterpart of Theorem~\ref{thmA} is used in the
proof of Theorem~\ref{thmA} and again in Section~\ref{sec:topology}.

\begin{lemma}\label{lem:punctured}
Let $n>r\ge 1$, and suppose that $r\le 2$ or $n\in\{r+1,r+2\}$. Let
$P_{n,r}=\VR(B^*_{n,r};r)$ and $G_{n,r}=H_{n,r}[B^*_{n,r}]$. Then
\[
  \im(G_{n,r})=\ly(P_{n,r})=\coll(P_{n,r})=\theta(P_{n,r})=2^r-1.
\]
\end{lemma}

\begin{proof}
By the construction in the proof of Proposition~\ref{prop:VR-link},
$G_{n,r}$ contains an induced matching with $2^r-1$ edges. Since
$P_{n,r}=\Ind(G_{n,r})$, Theorem~\ref{thm:ind-cochord} gives
\[
  2^r-1\le\im(G_{n,r})\le\ly(P_{n,r})\le\coll(P_{n,r})
  \le\theta(P_{n,r}).
\]
It remains to show that $\theta(P_{n,r})\le 2^r-1$.
\vskip 1mm
\noindent
\emph{Case $r=1$.} Here $B^*_{n,1}=\{e_1,\ldots,e_n\}$, and any two
of these vertices are at Hamming distance $2$. Hence $P_{n,1}$ is
zero-dimensional, and $\theta(P_{n,1})\le 1$ by
Lemma~\ref{lem:theta-basic}.
\vskip 1mm
\noindent
\emph{Case $r=2$.} Identify the vertices of $B^*_{n,2}$ with the
$1$-element and $2$-element subsets of $[n]$, and let $G=G_{n,2}$.
Two singletons are never adjacent in $G$; a singleton $\{i\}$ and a
$2$-element set $A$ are adjacent if and only if $i\notin A$; and two
$2$-element sets are adjacent if and only if they are disjoint.
Consequently, up to a permutation of coordinates, every edge $ab$ of
$G$ is of one of the following two types, where
$W=V(G)\setminus(N[a]\cup N[b])$:
\begin{enumerate}[label=\textup{(\roman*)}]
\item $a=\{1\}$ and $b=\{2,3\}$, in which case
$W=\{\{2\},\{3\},\{1,2\},\{1,3\}\}$;
\item $a=\{1,2\}$ and $b=\{3,4\}$ (possible only if $n\ge 4$), in
which case $$W=\{\{1,3\},\{1,4\},\{2,3\},\{2,4\}\}.$$
\end{enumerate}
In both cases $G[W]$ is a matching with two edges, so
$\alpha(G[W])=2$ and $\mu(G)=2$. Lemma~\ref{lem:common} gives
$\theta(P_{n,2})\le 3$.

\emph{Case $n=r+1$.} Here
$B^*_{r+1,r}=\F^{r+1}\setminus\{0,\one\}$. Hence $G_{r+1,r}$ is
obtained from the perfect matching $H_{r+1,r}$ by deleting both
endpoints of the edge $0\one$, so it is a perfect matching with
$2^r-1$ edges. By Lemma~\ref{lem:common}, $\theta(P_{r+1,r})=2^r-1$.

\emph{Case $n=r+2$.} Here
$B^*_{r+2,r}=\F^{r+2}\setminus(B(\one,1)\cup\{0\})$. Consider the
co-chordal cover $\{P_p:p\in\mathcal P\}$ of $H_{r+2,r}$ constructed
in the proof of Theorem~\ref{thm:VR-n-2}, and let
$p_0=\{0^{r+1},1^{r+1}\}$. Both vertices of $K_{p_0}$, namely $0$ and
$1^{r+1}0$, lie outside $B^*_{r+2,r}$, and so does the endpoint
$\one$ of the edge in $F_{p_0}$. Hence every edge of $P_{p_0}$ has
an endpoint outside $B^*_{r+2,r}$, and $P_{p_0}[B^*_{r+2,r}]$ is
edgeless. An induced subgraph of a co-chordal graph is co-chordal,
since its complement is an induced subgraph of a chordal graph.
Therefore the $2^r-1$ graphs $P_p[B^*_{r+2,r}]$, $p\ne p_0$, are
co-chordal, and they cover $E(G_{r+2,r})$. Thus
$\cochord(G_{r+2,r})\le 2^r-1$, and Theorem~\ref{thm:ind-cochord}
gives $\theta(P_{r+2,r})\le 2^r-1$.
\end{proof}

\subsection{Proof of Theorem~\texorpdfstring{\ref{thmA}}{A}}
\label{subsec:proofC}

\begin{proof}
Suppose that $r\le 2$ or $n\in\{r+1,r+2\}$. By
Lemma~\ref{lem:punctured} and Proposition~\ref{prop:VR-link}(a), we
have
$\ly(\VR(I_n;r))=\coll(\VR(I_n;r))=\theta(\VR(I_n;r))=2^r;$
for $n=r+2$, this also follows directly from
Theorem~\ref{thm:VR-n-2}. The equality for the regularity follows
from Proposition~\ref{prop:VR-flag}.
For Lew's invariant, the cases $r=1$ and $r=2$ are
Proposition~\ref{prop:VR-r1} and Theorem~\ref{thm:VR-r2},
respectively. If $n=r+1$, then $\VR(I_n;r)$ is the boundary of a
$2^r$-dimensional cross-polytope, so
\[
  2^r=\theta(\VR(I_n;r))\le d'(\VR(I_n;r))\le\dim\VR(I_n;r)+1=2^r
\]
by Theorem~\ref{thm:theta-dprime}.
Finally, the co-chordal cover number equals $2^r$ by
Proposition~\ref{prop:VR-r1} when $r=1$; because $H_{r+1,r}$ is a
perfect matching with $2^r$ edges when $n=r+1$; and by
Theorem~\ref{thm:VR-n-2} when $n=r+2$.
\end{proof}

\section{Upper bounds independent of \texorpdfstring{$n$}{n}}
\label{sec:uniform}

In this section we establish upper bounds for $\theta(\VR(I_n;r))$
that depend only on $r$, under Hypothesis~$(\mathrm S_r)$. We first
apply $(\mathrm S_r)$, the binary case of the skew Hamming set-pair
theorem of Gao and Zhao~\cite{GZ}, to Lew's invariant, and we then
improve the resulting bound by using links of faces. The escape
lemma (Lemma~\ref{lem:escape}) and the iterated link bound
(Lemma~\ref{lem:iterated-link}) do not depend on $(\mathrm S_r)$.

\subsection{Separated sequences and Lew's invariant}
\label{subsec:separated}

A sequence $v_1,\ldots,v_k$ of points of $\F^n$ is
\emph{$r$-separated} if there exist \emph{witnesses}
$w_1,\ldots,w_k\in\F^n$ such that
\begin{equation}\label{eq:skew}
  d_H(w_i,v_i)\ge r+1
  ~\text{and}~
  d_H(w_i,v_j)\le r~(1\le j<i\le k).
\end{equation}
These conditions imply that the points $v_i$ are distinct, and so are
the points $w_i$. However, the points $v_i$ need not be pairwise
$r$-close. Adding this condition yields a description of Lew's
invariant.

\begin{lemma}\label{lem:dprime-VR}
Let $n>r\ge 1$, and let $Z\subseteq\F^n$. Then $d'(\VR(Z;r))$ is the
maximum length of an $r$-separated sequence $v_1,\ldots,v_k$ with
witnesses $w_1,\ldots,w_k$, such that all $v_i$ and $w_i$ lie in $Z$
and $d_H(v_i,v_j)\le r$ for all $i$ and $j$. The empty sequence is
allowed.
\end{lemma}

\begin{proof}
Since $\VR(Z;r)=\Ind(H_{n,r}[Z])$, conditions~(i)-(iii) of
Lemma~\ref{lem:dprime-flag}, applied with $p_i=v_i$ and $q_i=w_i$,
are precisely the pairwise $r$-closeness of the points $v_i$ and the
inequalities~\eqref{eq:skew}.
\end{proof}

We use the following sharp skew Hamming set-pair inequality.

\begin{theorem}[{Gao and Zhao~\cite[Theorem~2.2]{GZ}}]\label{thm:GZ}
Let $A$ be an alphabet, and let $n>t\ge 0$. Suppose that
$(a_i,b_i)_{i=1}^m$ are pairs of words in $A^n$ such that
$d_H(a_i,b_i)\ge t+1$ for every $i$ and $d_H(a_i,b_j)\le t$ whenever
$i<j$. Then $m\le 2^{t+1}$. If $|A|\ge 2$, this bound is sharp.
\end{theorem}

Taking $A=\F$, $t=r$, $a_i=v_i$, and $b_i=w_i$, the
conditions~\eqref{eq:skew} become the hypotheses of
Theorem~\ref{thm:GZ}. Hence Theorem~\ref{thm:GZ} implies
Hypothesis~$(\mathrm S_r)$ for every $r\ge 1$. Since~\cite{GZ} is
currently a preprint, we state the results that depend on it as
conditional on $(\mathrm S_r)$.

\begin{theorem}\label{thm:VR-uniform}
Let $n>r\ge 1$, and assume Hypothesis~$(\mathrm S_r)$. Then
\[
  2^r\le\ly(\VR(I_n;r))\le\coll(\VR(I_n;r))\le\theta(\VR(I_n;r))
  \le d'(\VR(I_n;r))\le 2^{r+1}.
\]
Moreover, for every nonempty $Y\subseteq\F^n$, the complex
$\VR(Y;r)$ is homotopy equivalent to a finite CW complex of
dimension at most $2^{r+1}-1$.
\end{theorem}

\begin{proof}
By Lemma~\ref{lem:dprime-VR}, every Lew sequence of $\VR(I_n;r)$
gives an $r$-separated sequence whose points $v_i$ are pairwise
$r$-close. Hence $(\mathrm S_r)$ gives
$d'(\VR(I_n;r))\le 2^{r+1}$, and the remaining inequalities follow
from Theorem~\ref{thm:VR-matching} and
Theorem~\ref{thmB}. For $Y\subseteq\F^n$,~\ref{F:mono}
gives $\theta(\VR(Y;r))\le\theta(\VR(I_n;r))\le 2^{r+1}$, and the
homotopy assertion follows from Lemma~\ref{lem:theta-homotopy}.
\end{proof}

The bound $2^{r+1}$ is sharp for $r$-separated sequences. Indeed, let
$Q$ be an $(r+1)$-dimensional coordinate subcube, list its vertices
as $v_1,\ldots,v_{2^{r+1}}$, and let $w_i$ be the antipode of $v_i$
in $Q$. Then $d_H(w_i,v_i)=r+1$, and $d_H(w_i,v_j)\le r$ for all
$i\ne j$. However, the points $v_i$ are not pairwise $r$-close, so
this construction does not show that the bound is sharp for $d'$.
Whether the additional closeness condition reduces the bound to
$2^r$ is the subject of Conjecture~\ref{conj:skew}.

\subsection{The escape and padding lemmas}
\label{subsec:escape}

We improve the uniform bound by extending separated sequences in
links of faces. The key ingredient is the following lemma.

\begin{lemma}[escape lemma]\label{lem:escape}
Let $n>r\ge 1$, let $u\in\F^n$, and let
$A\subseteq\F^n\setminus\{u\}$ satisfy $|A|\le 3$ and
$\diam(A\cup\{u\})\le r$. Then there exists $z\in\F^n$ such that
$d_H(z,u)\ge r+1$ and $d_H(z,a)\le r$ for every $a\in A$.
\end{lemma}

\begin{proof}
By translation, we may assume that $u=0$, and we identify each
vector with its support. If $A=\emptyset$, let $z$ be any
$(r+1)$-element subset of $[n]$. Otherwise, write
$A=\{S_1,S_2,S_3\}$, allowing repetitions. Then $1\le|S_i|\le r$ and
$|S_i\triangle S_j|\le r$ for all $i$ and $j$. Let
$U=S_1\cup S_2\cup S_3$.
If $|U|\le r+1$, extend $U$ to a set $z$ with $|z|=r+1$. Then
$S_i\subseteq z$, so $d_H(z,S_i)=r+1-|S_i|\le r$ for each $i$.
Suppose now that $|U|\ge r+2$. For $\{i,j,k\}=\{1,2,3\}$, let
\[
  P_i=S_i\setminus(S_j\cup S_k)
  ~\text{and}~
  Q_{jk}=(S_j\cap S_k)\setminus S_i,
\]
and write $p_i=|P_i|$, $q_{jk}=|Q_{jk}|$, and
$c_i=|U\setminus S_i|=p_j+p_k+q_{jk}$. If $c_i\le r$ for every $i$,
then $z=U$ has the required properties. Otherwise, after relabeling
we may assume that $c_1=\max_i c_i\ge r+1$. Set $\varepsilon_i=c_i-r$, so that
$\varepsilon_1\ge 1$ and $\varepsilon_2,\varepsilon_3 \le \varepsilon_1$.

Since $|S_1\triangle S_2|=p_1+p_2+q_{13}+q_{23}\le r\le c_1-1$, we
obtain $p_1+q_{13}\le p_3-1$; similarly,
$|S_1\triangle S_3|\le r$ gives $p_1+q_{12}\le p_2-1$. Consequently,
\[
  \varepsilon_2+\varepsilon_3=2p_1+p_2+p_3+q_{12}+q_{13}-2r
  \le 2(p_2+p_3)-2-2r\le -2,
\]
where the last inequality holds because
$p_2+p_3\le|S_2\triangle S_3|\le r$. Moreover,
$\varepsilon_1=p_2+p_3+q_{23}-r\le\min\{p_2,p_3\}$, since
$p_3+q_{23}\le|S_3|\le r$ and $p_2+q_{23}\le|S_2|\le r$.
Consider the integer interval
$I=[\max\{\varepsilon_2,-\varepsilon_1\},\,\min\{-\varepsilon_3,\varepsilon_1\}]$. It is nonempty, because
$\varepsilon_2\le -\varepsilon_3$, $\varepsilon_2\le \varepsilon_1$, $\varepsilon_3\le \varepsilon_1$, and $\varepsilon_1\ge 1$. We claim
that $I$ contains an integer $d$ with $d\equiv \varepsilon_1\pmod 2$. This is
clear if $I$ contains at least two integers. If $I=\{m\}$, then
$m\in\{\varepsilon_1,-\varepsilon_1\}$, since otherwise $m=\varepsilon_2=-\varepsilon_3$, contradicting
$\varepsilon_2+\varepsilon_3\le -2$.
Set $y_2=(\varepsilon_1-d)/2$ and $y_3=(\varepsilon_1+d)/2$. These are nonnegative
integers with $y_2+y_3=\varepsilon_1\le\min\{p_2,p_3\}$. Choose
$D\subseteq P_2\cup P_3$ with $|D\cap P_2|=y_2$ and
$|D\cap P_3|=y_3$, and let $z=U\setminus D$. Since $D$ is disjoint
from $S_1$ and $|U|=|S_1|+c_1$, we have $S_1\subseteq z$ and
\[
  |z|=|U|-\varepsilon_1=|S_1|+r\ge r+1,
  ~
  d_H(z,S_1)=|z|-|S_1|=r.
\]
Moreover, $z\setminus S_2=(U\setminus S_2)\setminus(D\cap P_3)$ and
$S_2\setminus z=D\cap P_2$, so
\[
  d_H(z,S_2)=c_2-y_3+y_2=c_2-d\le r,
  ~
  d_H(z,S_3)=c_3-y_2+y_3=c_3+d\le r,
\]
where the inequalities follow from $\varepsilon_2\le d\le -\varepsilon_3$. Thus $z$ has
the required properties.
\end{proof}

The restriction $|A|\le 3$ cannot be removed in general. For
example, let $r=2$, $n\ge 4$, $u=0$, and $A=\{e_1,e_2,e_3,e_4\}$. If
$|z|\ge 3$, then $d_H(z,e_i)=|z|-1$ when $i\in z$ and
$d_H(z,e_i)=|z|+1$ otherwise; hence $d_H(z,e_i)\le 2$ for all
$i\le 4$ would force $\{1,2,3,4\}\subseteq z$ and $|z|\le 3$, which
is impossible.

\begin{lemma}\label{lem:padding}
Let $n>r\ge 1$, let $\Delta=\VR(I_n;r)$, and let
$\sigma=\{u_1,\ldots,u_s\}$ be a nonempty face of $\Delta$. Suppose
that there exist points $y_1,\ldots,y_s,x_1,\ldots,x_s\in\F^n$ such
that, for $1\le\ell \le s$,
\begin{equation}\label{eq:padding}
\begin{aligned}
  d_H(y_{\ell},u_{\ell})&\ge r+1, &~ d_H(y_{\ell},u_j)&\le r~(j<\ell),\\
  d_H(x_{\ell},u_{\ell})&\ge r+1, &~ d_H(x_{\ell},u_j)&\le r~(j>\ell).
\end{aligned}
\end{equation}
If Hypothesis~$(\mathrm S_r)$ holds, then
$d'(\lk(\Delta,\sigma))\le 2^{r+1}-2s$. Moreover, independently of
$(\mathrm S_r)$, such points exist for every face with $s\le 4$ and
every ordering of its vertices.
\end{lemma}

\begin{proof}
Since $\Delta$ is flag, $\lk(\Delta,\sigma)=\VR(Z;r)$, where $Z$ is
the set of points outside $\sigma$ that are $r$-close to every vertex
of $\sigma$. Let $k=d'(\lk(\Delta,\sigma))$. By
Lemma~\ref{lem:dprime-VR}, there is an $r$-separated sequence
$v_1,\ldots,v_k$ of pairwise $r$-close points of $Z$, with witnesses
$w_1,\ldots,w_k\in Z$. Consider the sequence
\[
  u_1,\ldots,u_s,\ v_1,\ldots,v_k,\ x_1,\ldots,x_s,
\]
with witnesses
\[
  y_1,\ldots,y_s,\ w_1,\ldots,w_k,\ u_1,\ldots,u_s,
\]
respectively. We claim that this sequence is $r$-separated.

The conditions on each point and its own witness follow
from~\eqref{eq:padding} and~\eqref{eq:skew}. It remains to check
that each witness is $r$-close to all preceding points. The witness
$y_{\ell}$ is $r$-close to the preceding points $u_j$, $j<\ell,$
by~\eqref{eq:padding}. The witness $w_i$ is $r$-close to every
vertex of $\sigma$, since $w_i\in Z$, and to the preceding points
$v_j$ by~\eqref{eq:skew}. Finally, the witness $u_{\ell}$ of $x_{\ell}$ is
$r$-close to every vertex of $\sigma$, since $\sigma$ is a face, to
every $v_i$, since $v_i\in Z$, and to the preceding points $x_j$,
$j<\ell$, by~\eqref{eq:padding}. This proves the claim.

The points of this sequence need not be pairwise $r$-close, so here
the full strength of Hypothesis~$(\mathrm S_r)$ is used. It gives
$k+2s\le 2^{r+1}$, as required.
Now suppose that $s\le 4$, and fix $\ell$. The sets
$\{u_j:j<\ell\}$ and $\{u_j:j>\ell\}$ each have at most three elements,
and each of them, together with $u_{\ell}$, has diameter at most $r$,
since $\sigma$ is a face. Applying Lemma~\ref{lem:escape} to these
two sets yields the required points $y_{\ell}$ and $x_{\ell}$.
\end{proof}

\subsection{The improved bound}
\label{subsec:padded}

To pass from links of faces to the whole complex, we iterate the
bound of Lemma~\ref{lem:link}.

\begin{lemma}[iterated link bound]\label{lem:iterated-link}
Let $\X$ be a simplicial complex, and let $s\ge 1$. Then
$\theta(\X)$ is at most the maximum of the following numbers:
\begin{enumerate}[label=\textup{(\roman*)}]
\item $|\tau|$, for faces $\tau\in\X$ \textup{(}including
$\tau=\emptyset$\textup{)} with $|\tau|<s$ such that
$\lk(\X,\tau)$ is a simplex;
\item $s+\theta(\lk(\X,\sigma))$, for faces $\sigma\in\X$ with
$|\sigma|=s$.
\end{enumerate}
\end{lemma}

\begin{proof}
We proceed by induction on $s$. For $s=1$, if $\X$ is a simplex,
then the face $\tau=\emptyset$ contributes $0=\theta(\X)$ to~(i);
otherwise, the claim is Lemma~\ref{lem:link}.
Suppose that the assertion holds for some $s\ge 1$. The numbers in~(i)
for $s$ also occur in~(i) for $s+1$. Let $\sigma$ be a face with
$|\sigma|=s$. If $\lk(\X,\sigma)$ is a simplex, then
$s+\theta(\lk(\X,\sigma))=s=|\sigma|$, which occurs in~(i) for
$s+1$. Otherwise, Lemma~\ref{lem:link} and the identity
$\lk(\lk(\X,\sigma),v)=\lk(\X,\sigma\cup\{v\})$ give
\[
  s+\theta(\lk(\X,\sigma))
  \le\max_{v}\bigl\{s+1+\theta(\lk(\X,\sigma\cup\{v\}))\bigr\},
\]
where $v$ ranges over the vertices of $\lk(\X,\sigma)$. These
numbers occur in~(ii) for $s+1$, which completes the induction.
\end{proof}

We now obtain the improved uniform bound.

\begin{theorem}\label{thm:VR-padded}
Let $n>r\ge 2$, and assume Hypothesis~$(\mathrm S_r)$. Then, for
every nonempty $Y\subseteq\F^n$,
\[
  \ly(\VR(Y;r))\le\coll(\VR(Y;r))\le\theta(\VR(Y;r))
  \le\theta(\VR(I_n;r))\le 2^{r+1}-4.
\]
Moreover, $\VR(Y;r)$ is homotopy equivalent to a finite CW complex
of dimension at most $2^{r+1}-5$.
\end{theorem}

\begin{proof}
Let $\Delta=\VR(I_n;r)$. For every face $\sigma$ of $\Delta$ with
$|\sigma|=4$, Lemma~\ref{lem:padding} gives
$d'(\lk(\Delta,\sigma))\le 2^{r+1}-8$, and hence
$\theta(\lk(\Delta,\sigma))\le 2^{r+1}-8$ by
Theorem~\ref{thm:theta-dprime}. Applying
Lemma~\ref{lem:iterated-link} with $s=4$, we obtain
$\theta(\Delta)\le\max\{3,\,2^{r+1}-4\}=2^{r+1}-4,$
where the equality holds because $r\ge 2$. The remaining
inequalities follow from~\ref{F:leray} and~\ref{F:mono}, and the
homotopy assertion follows from Lemma~\ref{lem:theta-homotopy}.
\end{proof}

\subsection{Explicit bounds}
\label{subsec:explicit}

Recall that $\kappa(n,r)=\dim\VR(I_n;r)+1$ is given by Kleitman's
diametric theorem~\cite{Kle}, and that
$U(n,r)=\min\{\kappa(n,r),2^{r+1}-4\}$ for $n>r\ge 2$.

\begin{corollary}\label{cor:best-bound}
Let $n>r\ge 1$.
\begin{enumerate}[label=\textup{(\alph*)}]
\item Unconditionally,
$\coll(\VR(I_n;r))\le\theta(\VR(I_n;r))\le d'(\VR(I_n;r))
\le\kappa(n,r)$.
\item If $r\ge 2$ and Hypothesis~$(\mathrm S_r)$ holds, then
$\theta(\VR(I_n;r))\le U(n,r)$ and
$d'(\VR(I_n;r))\le\min\{\kappa(n,r),2^{r+1}\}$.
\end{enumerate}
For $r\ge 2$, let $n_0(r)$ be the largest integer $n>r$ such that
$\kappa(n,r)\le 2^{r+1}-4$. Then
\[
  U(n,r)=
  \begin{cases}
    \kappa(n,r) & \text{if } r<n\le n_0(r),\\
    2^{r+1}-4 & \text{if } n>n_0(r).
  \end{cases}
\]
For example, $n_0(4)=6$, since $\kappa(6,4)=22$ and
$\kappa(7,4)=29$. Consequently, the bound $\theta(\VR(I_n;r))\le U(n,r)$
holds unconditionally whenever $r=2$ or $n\le n_0(r)$.
\end{corollary}

\begin{proof}
By Theorem~\ref{thm:theta-dprime} and Lemma~\ref{lem:dprime-dim},
$\theta(\VR(I_n;r))\le d'(\VR(I_n;r))\le\kappa(n,r)$, which
proves~(a). Part~(b) follows from~(a) together with
Theorems~\ref{thm:VR-uniform} and~\ref{thm:VR-padded}.
Fix $r\ge 2$. By Kleitman's formula, $\kappa(n,r)$ is strictly
increasing and unbounded as a function of $n$. Since
$\kappa(r+1,r)=2^r\le 2^{r+1}-4$, the integer $n_0(r)$ exists, and
the stated description of $U(n,r)$ follows.  
Finally, if $r=2$, then $\theta(\VR(I_n;2))\le 4=U(n,2)$ by
Theorem~\ref{thm:VR-r2}, and if $r<n\le n_0(r)$, then
$U(n,r)=\kappa(n,r)$ and the bound $\theta(\VR(I_n;r))\le U(n,r)$
follows from~(a).
\end{proof}

For fixed $r$, the quantity $\kappa(n,r)$ grows on the order of
$n^{\lfloor r/2\rfloor}$. Bendersky, Elia, and
Grbi\'c~\cite[Theorem~4.2]{BEG} obtained coconnectivity bounds for
the cohomology of the full complex $\VR(I_n;r)$; these do not
directly control the Leray number or the collapsibility number,
which depend on all induced subcomplexes. By contrast, under
Hypothesis~$(\mathrm S_r)$, Corollary~\ref{cor:best-bound} gives
\[
  2^r\le\coll(\VR(I_n;r))\le\theta(\VR(I_n;r))\le 2^{r+1}-4
  ~(n>r\ge 2),
\]
uniformly in $n$. The smallest case not covered by
Theorem~\ref{thmA} or by Shukla's results is $(n,r)=(7,4)$, for
which we have unconditionally
\[
  16\le\coll(\VR(I_7;4))\le\theta(\VR(I_7;4))\le 29=\kappa(7,4),
\]
and $\theta(\VR(I_7;4))\le 28$ under $(\mathrm S_4)$.

\subsection{Co-chordal covers for general \texorpdfstring{$n$}{n}}
\label{subsec:codes-short}

By Theorem~\ref{thmB}, Conjecture~\ref{conj:theta} for a pair $(n,r)$
would follow from a cover of $H_{n,r}$ by $2^r$ co-chordal subgraphs,
as constructed for $n=r+2$ in Theorem~\ref{thm:VR-n-2}. For general
$n$, Lemma~\ref{lem:star-clique} with $F=\emptyset$ shows that, for
every clique $K$ of a graph $G$, the edges of $G$ meeting $K$ form a
co-chordal subgraph (compare~\cite[Theorem~2]{W14}). If $C\le\F^n$ is a binary linear code of
dimension $k$ and minimum distance at least $r+1$, its $2^{n-k}$
cosets partition $\F^n$ into cliques of $H_{n,r}$, so
$\cochord(H_{n,r})\le 2^{n-k}$. Shortened (extended) BCH
codes~\cite[Chapter~7, Section~6]{MacS} then give
$\cochord(H_{n,r})\le 2(2n)^{\lfloor r/2\rfloor}$. This bound is
polynomial in $n$ for fixed $r$, but it is at least $\kappa(n,r)$
(indeed, $\sum_{i=0}^{t}\binom{m}{i}\le(m+1)^t$ for all integers
$m,t\ge 0$, so Kleitman's formula gives
$\kappa(n,r)\le 2(n+1)^{\lfloor r/2\rfloor}\le 2(2n)^{\lfloor r/2\rfloor}$), so
it does not improve the bound $\theta(\VR(I_n;r))\le\kappa(n,r)$.
We do not know whether $\cochord(H_{n,r})=2^r$ for any pair with
$r\ge 2$ and $n\ge r+3$; the smallest such case is $(n,r)=(5,2)$.

\subsection{Proof of Theorem~\texorpdfstring{\ref{thmD}}{D}}
\label{subsec:proofB}

\begin{proof}
(a) By Theorem~\ref{thm:VR-matching}, an $(r+1)$-dimensional
coordinate subcube induces a perfect matching with $2^r$ edges, so
$\im(H_{n,r})\ge 2^r$. The reverse inequality is
Corollary~\ref{thm:im-exact}, which uses only the symmetric theorem
of Alon, Jin, and Sudakov (Theorem~\ref{thm:AJS}). The inequality
$\im(H_{n,r})\le\ly_\Bbbk(\VR(I_n;r))$ is part of
Theorem~\ref{thm:ind-cochord}, and $\ly_\Bbbk\le\coll\le\theta$
is~\ref{F:leray}. Theorem~\ref{thm:theta-dprime} and
Lemma~\ref{lem:dprime-dim} give
$\theta(\VR(I_n;r))\le d'(\VR(I_n;r))\le\kappa(n,r)$. None of
these steps uses Theorem~\ref{thm:GZ}.

(b) Assume $(\mathrm S_r)$. By Theorem~\ref{thm:VR-uniform} and
part~(a), $d'(\VR(I_n;r))\le\min\{\kappa(n,r),2^{r+1}\}$. By
Theorem~\ref{thm:VR-padded}, $\theta(\VR(I_n;r))\le 2^{r+1}-4$;
combined with~(a), this gives $\theta(\VR(I_n;r))\le U(n,r)$, as
recorded in Corollary~\ref{cor:best-bound}(b). The homotopy
assertion is part of Theorem~\ref{thm:VR-padded}.
\end{proof}

\subsection{A close-skew set-pair conjecture}
\label{subsec:close-skew}

Lemma~\ref{lem:dprime-VR} suggests strengthening the Gao--Zhao bound
under the additional assumption that the first points of the
ordered pairs are pairwise $r$-close.

\begin{conjecture}\label{conj:skew}
Let $n>r\ge 1$, and let $v_1,\ldots,v_k\in\F^n$ satisfy
$d_H(v_i,v_j)\le r$ for all $i$ and $j$. If there exist
$w_1,\ldots,w_k\in\F^n$ such that $d_H(w_i,v_i)\ge r+1$ for every $i$
and $d_H(w_i,v_j)\le r$ whenever $j<i$, then $k\le 2^r$.
\end{conjecture}

The conjectured bound would be sharp. Indeed, let $Q$ be an
$r$-dimensional coordinate subcube contained in an
$(r+1)$-dimensional coordinate subcube $R$, let $v_1,\ldots,v_{2^r}$
be the vertices of $Q$, and let $w_i$ be the antipode of $v_i$ in
$R$. The points $v_i$ are pairwise $r$-close, $d_H(w_i,v_i)=r+1$, and
$d_H(w_i,v_j)=r+1-d_H(v_i,v_j)\le r$ for $i\ne j$.
By Lemma~\ref{lem:dprime-VR}, Conjecture~\ref{conj:skew} for a given
pair $(n,r)$ is equivalent to the inequality
$d'(\VR(I_n;r))\le 2^r$. Since
$d'(\VR(I_n;r))\ge\theta(\VR(I_n;r))\ge 2^r$ by
Theorems~\ref{thm:theta-dprime} and~\ref{thm:VR-matching}, it is
therefore equivalent to
\begin{equation}\label{eq:skew-dprime}
  d'(\VR(I_n;r))=2^r.
\end{equation}
Together with Theorem~\ref{thmB}, this would imply
$\ly(\VR(I_n;r))=\coll(\VR(I_n;r))=\theta(\VR(I_n;r))=2^r$, and hence
Conjectures~\ref{conj:theta} and~\ref{conj:shukla}. By the exact
values established in Section~\ref{sec:small},
Conjecture~\ref{conj:skew} holds for $r\le 2$ and for $n=r+1$.

It is useful to distinguish three Hamming set-pair problems. In the
\emph{symmetric} problem, the cross-distance condition is required
for all $i\ne j$; Alon, Jin, and Sudakov~\cite{AJS} proved the sharp
bound $2^{r+1}$, which yields $\im(H_{n,r})=2^r$. In the
\emph{skew} problem, posed by Alon, Jin, and
Sudakov~\cite[Section~5]{AJS}, the condition is required only for
$j<i$; Gao and Zhao~\cite{GZ} proved, in a preprint, the same sharp
bound $2^{r+1}$. Both bounds are attained by pairing the vertices of
an $(r+1)$-dimensional coordinate subcube with their antipodes. The
\emph{close-skew} problem additionally requires the points $v_i$ to
be pairwise $r$-close, and Conjecture~\ref{conj:skew} predicts the
smaller bound $2^r$.
The padding lemma provides further evidence. Under
$(\mathrm S_r)$, if $\sigma$ is a face of $\VR(I_n;r)$ with
$|\sigma|\le 4$, then every close-skew system whose points and
witnesses lie in the vertex set of $\lk(\VR(I_n;r),\sigma)$ has at
most $2^{r+1}-2|\sigma|$ pairs.

\section{Topology and commutative algebra}
\label{sec:topology}

In this section we derive consequences of the preceding bounds for
the topology of the complexes $\VR(Y;r)$ with $Y\subseteq\F^n$ and
for their Stanley--Reisner ideals, and we prove Theorem~\ref{thmC}.

\subsection{Homotopy dimension and top homology}
\label{subsec:top-homology}

Every $(r+1)$-dimensional coordinate subcube supports a sphere of
dimension $2^r-1$. Coordinate projections give the following
retraction property, in the spirit of Adams and
Virk~\cite[Section~3]{AV}.

\begin{proposition}\label{prop:retraction}
Let $n>r\ge 1$, and let $Q$ be an $(r+1)$-dimensional coordinate
subcube with $Q\subseteq Y\subseteq\F^n$. Then $\VR(Q;r)$ is
homeomorphic to $S^{2^r-1}$ and is a simplicial retract of
$\VR(Y;r)$. Consequently,
$\widetilde H_{2^r-1}(\VR(Y;r);\mathbb Z)$ contains $\mathbb Z$ as a
direct summand.
\end{proposition}

\begin{proof}
Let $J\subseteq[n]$ be the set of coordinates that are fixed on $Q$,
with values $c_j$ for $j\in J$. Define $\pi:\F^n\to Q$ by replacing
the $j$th coordinate of each vector by $c_j$ for every $j\in J$. Then
$d_H(\pi(\alpha),\pi(\beta))\le d_H(\alpha,\beta)$ for all
$\alpha,\beta\in\F^n$, so $\pi$ maps faces of $\VR(Y;r)$ to faces of
$\VR(Q;r)$ and induces a simplicial map
$\pi_*:\VR(Y;r)\to\VR(Q;r)$. Since $\pi$ fixes $Q$ pointwise,
$\pi_*$ is a retraction. By Theorem~\ref{thm:VR-matching},
$\VR(Q;r)$ is the boundary of a $2^r$-dimensional cross-polytope,
and hence is homeomorphic to $S^{2^r-1}$. The retraction splits the
map on reduced integral homology induced by the inclusion, which
proves the last assertion.
\end{proof}

\begin{theorem}\label{thm:topology}
Let $n>r\ge 1$, and let $\emptyset\ne Y\subseteq\F^n$. Let
$U^\star=\kappa(n,r)$ in general, and $U^\star=U(n,r)$ if $r\ge 2$
and Hypothesis~$(\mathrm S_r)$ holds, where $\kappa$ and $U$ are as
in Section~\ref{sec:intro}.
\begin{enumerate}[label=\textup{(\alph*)}]
\item The complex $\VR(Y;r)$ is homotopy equivalent to a finite CW
complex of dimension at most $U^\star-1$; under $(\mathrm S_r)$ and
for $r\ge 2$, this dimension is at most $2^{r+1}-5$. Consequently,
$\widetilde H_i(\VR(Y;r);\mathbb Z)=0$ for all $i\ge U^\star$, and,
over every field, $\reg(S_Y/I_{\VR(Y;r)})\le U^\star$, where
$S_Y=\Bbbk[x_\alpha:\alpha\in Y]$. If $I_{\VR(Y;r)}\ne 0$, then
$\reg I_{\VR(Y;r)}\le U^\star+1$.
\item If $r\le 2$ or $n\le r+2$, then $\VR(Y;r)$ is homotopy
equivalent to a finite CW complex of dimension at most $2^r-1$.
\item If $r\le 3$ or $n\le r+2$, then
$\widetilde H_{2^r-1}(\VR(I_n;r);\mathbb Z)$ is a nonzero free
abelian group, and $\widetilde H_i(\VR(I_n;r);\mathbb Z)=0$ for all
$i\ge 2^r$. Moreover, $\reg I_{\VR(I_n;r)}=2^r+1$ over every field,
and the coconnectivity of $\VR(I_n;r)$ in the sense
of~\cite[Definition~4.1]{BEG}, that is, the largest integer $d$ such
that $H^{d-1}(\VR(I_n;r);\mathbb Z)\ne 0$, equals $2^r$.
\end{enumerate}
\end{theorem}

\begin{proof}
(a) By~\ref{F:mono} and Corollary~\ref{cor:best-bound},
$\theta(\VR(Y;r))\le\theta(\VR(I_n;r))\le U^\star$, and
Lemma~\ref{lem:theta-homotopy} gives the bound on the CW dimension
and the vanishing of integral homology. Moreover,
$\ly(\VR(Y;r))\le U^\star$ over every field by~\ref{F:leray}, and the
assertions about regularity follow from Hochster's
formula~\eqref{eq:hochster-betti}.

(b) By Theorem~\ref{thmA}, $\theta(\VR(I_n;r))=2^r$ in the stated
ranges. The conclusion follows from~\ref{F:mono} and
Lemma~\ref{lem:theta-homotopy}.

(c) Let $\Delta=\VR(I_n;r)$ and $d=2^r$. By Theorem~\ref{thmA}, or by
Shukla's result~\cite{Shu} when $r=3$, we have $\coll(\Delta)=d$.
By Wegner's theorem~\cite{Weg}, $\Delta$ is homotopy equivalent to a
finite CW complex of dimension at most $d-1$. Hence
$\widetilde H_i(\Delta;\mathbb Z)=0$ for $i\ge d$, and
$\widetilde H_{d-1}(\Delta;\mathbb Z)$ is free abelian, being a
subgroup of the top cellular chain group; it is nonzero by
Proposition~\ref{prop:retraction}. Since $d-1\ge 1$, the universal
coefficient theorem gives $H^{d-1}(\Delta;\mathbb Z)\ne 0$, while
$H^i(\Delta;\mathbb Z)=0$ for $i\ge d$ by the bound on the CW
dimension. Thus the coconnectivity of $\Delta$ equals $d$. Finally,
Theorem~\ref{thm:VR-matching} and~\ref{F:leray} give
$d\le\ly(\Delta)\le\coll(\Delta)=d$ over every field. Since
$I_\Delta\ne 0$, Hochster's formula yields $\reg I_\Delta=d+1$.
\end{proof}

Shukla~\cite{Shu} conjectured that, for $r\ge 2$ and $n\ge r+2$, the
reduced integral homology of $\VR(I_n;r)$ is nonzero exactly in
degrees $r+1$ and $2^r-1$. The prediction for the lower degree fails
in general by work of Bendersky and Grbi\'c~\cite{BG}. Part~(c) shows
that the prediction for the top degree holds whenever $r\le 3$ or
$n\le r+2$; for $n=r+2$, this closes the gap between the upper bounds
of~\cite[Theorem~4.2]{BEG} and the nonvanishing in degree $2^r-1$
recorded in~\cite[Table~2]{BEG}.

In the parameter ranges $r\le 2$ or $n\le r+2$, we can show further
that homology in the top degree cannot disappear when vertices are
added.

\begin{theorem}\label{thm:nonvanishing}
Let $n>r\ge 1$, suppose that $r\le 2$ or $n\le r+2$, and let
$m=2^r-1$.
\begin{enumerate}[label=\textup{(\alph*)}]
\item For all $Y'\subseteq Y\subseteq\F^n$, the inclusion
$\VR(Y';r)\hookrightarrow\VR(Y;r)$ induces a split injection
\[
  \widetilde H_m(\VR(Y';r);\mathbb Z)\longrightarrow
  \widetilde H_m(\VR(Y;r);\mathbb Z)
\]
with free abelian cokernel.
\item If $Y$ contains a subset $Y'$ such that $\VR(Y';r)$ is
homeomorphic to $S^m$, for instance an $(r+1)$-dimensional
coordinate subcube, then
$\widetilde H_m(\VR(Y;r);\mathbb Z)\ne 0$, and $m$ is the largest
degree in which $\VR(Y;r)$ has nonzero reduced homology.
\end{enumerate}
\end{theorem}

\begin{proof}
(a) If $f$ and $g$ are composable injections of finitely generated
abelian groups with free cokernels, then $gf$ is injective and there
is an exact sequence
$0\to\operatorname{coker}f\to\operatorname{coker}(gf)
\to\operatorname{coker}g\to 0$, so $\operatorname{coker}(gf)$ is
free and $gf$ splits. Hence it suffices to treat the case
$Y=Y'\cup\{y\}$ with $y\notin Y'$.
Let $\X=\VR(I_n;r)$, $A=\X[Y']$, $B=\X[Y]$, and
$\Lambda=\lk(B,y)$. Then $B=A\cup(y*\Lambda)$ and
$A\cap(y*\Lambda)=\Lambda$. Since the cone $y*\Lambda$ is acyclic,
the reduced Mayer--Vietoris sequence contains the exact segment
\[
  \widetilde H_m(\Lambda;\mathbb Z)\longrightarrow
  \widetilde H_m(A;\mathbb Z)\longrightarrow
  \widetilde H_m(B;\mathbb Z)\longrightarrow
  \widetilde H_{m-1}(\Lambda;\mathbb Z).
\]
The link $\Lambda$ is an induced subcomplex of $\lk(\X,y)$, which is
isomorphic to $P_{n,r}$. By Lemma~\ref{lem:punctured}
and~\ref{F:mono}, $\theta(\Lambda)\le m$, so by
Lemma~\ref{lem:theta-homotopy}, $\Lambda$ is homotopy equivalent to
a finite CW complex of dimension at most $m-1$. Consequently,
$\widetilde H_m(\Lambda;\mathbb Z)=0$ and
$\widetilde H_{m-1}(\Lambda;\mathbb Z)$ is free abelian. (If
$\Lambda=\{\emptyset\}$, both groups vanish because $m\ge 1$.) By
exactness, the map $\widetilde H_m(A;\mathbb Z)\to
\widetilde H_m(B;\mathbb Z)$ is injective, and its cokernel embeds
in $\widetilde H_{m-1}(\Lambda;\mathbb Z)$. The cokernel is
therefore free abelian, and the injection splits.

(b) By~(a), the fundamental class of the sphere $\VR(Y';r)$ maps to
a nonzero class in $\widetilde H_m(\VR(Y;r);\mathbb Z)$.
Theorem~\ref{thm:topology}(b) gives the vanishing of homology above
degree $m$.
\end{proof}

Consequently, if $r\le 2$ or $n\le r+2$ and $Y$ contains an
$(r+1)$-dimensional coordinate subcube, then both the homotopy
dimension of $\VR(Y;r)$ and the largest degree of its nonzero
homology equal $2^r-1$.
Such spheres need not arise from coordinate subcubes. For example,
let $(n,r)=(4,2)$, and let $Y$ be the set of the eight vectors of
even weight, which is not a coordinate subcube. Two
distinct elements of $Y$ are at distance greater than $2$ if and only
if they are antipodal. Thus $H_{4,2}[Y]$ is a perfect matching with
four edges, and $\VR(Y;2)$ is a $3$-sphere. For the full hypercube,
Adams and Virk~\cite[Theorem~4.1]{AV} proved that
\[
  \rank\widetilde H_{2^r-1}(\VR(I_n;r);\mathbb Z)
  \ge 2^{n-r-1}\binom{n}{r+1}
  ~(n>r\ge 2).
\]

\subsection{Vertex decomposability and depth}
\label{subsec:VD}

Recall that a simplicial complex $\X$ is \emph{vertex decomposable}
if either $\X$ is a simplex, or there is a vertex $v$ of $\X$ such
that $\X\setminus v$ and $\lk(\X,v)$ are vertex decomposable and no
face of $\lk(\X,v)$ is a facet of $\X\setminus v$; such a vertex is
called a \emph{shedding vertex}~\cite{BW,PB}. Vertex decomposable
complexes are shellable, and hence sequentially Cohen--Macaulay.
B{\i}y{\i}ko\u{g}lu and Civan~\cite[Theorem~30]{BC} proved that
\begin{equation}\label{eq:BC-vd}
  \theta(\X)=\coll(\X)=\ly(\X)
\end{equation}
for every vertex decomposable complex $\X$.

We also recall the algebraic notions used below; see~\cite{BH}. Let
$S=\Bbbk[x_v:v\in V]$ be standard graded, let $\mathfrak m$ be its
homogeneous maximal ideal, and let $M\ne 0$ be a finitely generated
graded $S$-module. The \emph{depth} of $M$, denoted $\depth M$, is
the common length of all maximal $M$-regular sequences contained in
$\mathfrak m$; equivalently,
$\depth M=\min\{i:H^i_{\mathfrak m}(M)\ne 0\}$, where
$H^i_{\mathfrak m}(M)$ denotes local cohomology with support in
$\mathfrak m$. The \emph{projective dimension} $\pd_S M$ is the length
of a minimal graded free resolution of $M$ over $S$, that is,
$\pd_S M=\max\{i:\beta_{i,j}(M)\ne 0\text{ for some }j\}$. These two
invariants are related by the Auslander--Buchsbaum formula
\begin{equation*}\label{eq:AB}
  \pd_S M+\depth M=|V|.
\end{equation*}
One always has $\depth M\le\dim M$, where $\dim M$ denotes the Krull
dimension, and $M$ is \emph{Cohen--Macaulay} if equality holds. A
simplicial complex $\X$ on $V$ is \emph{Cohen--Macaulay} over $\Bbbk$
if the Stanley--Reisner ring $S/I_\X$ is Cohen--Macaulay. Recall that
$\dim(S/I_\X)=\dim\X+1$, and that, by Hochster's formula for local
cohomology~\cite[Theorem~5.3.8]{BH},
\begin{equation*}\label{eq:hochster-depth}
  \depth(S/I_\X)=\min\bigl\{|\sigma|+j+1:\ \sigma\in\X,\
  \widetilde H_j(\lk(\X,\sigma);\Bbbk)\ne 0\bigr\}.
\end{equation*}
In particular, by Reisner's criterion~\cite[Corollary~5.3.9]{BH},
$\X$ is Cohen--Macaulay over $\Bbbk$ if and only if
$\widetilde H_i(\lk(\X,\sigma);\Bbbk)=0$ for every face $\sigma$ of
$\X$ and every $i<\dim\lk(\X,\sigma)$; Cohen--Macaulay complexes are
pure.

We now characterize vertex decomposability, the Cohen--Macaulay
property, and purity for the complexes $\VR(I_n;r)$.

\begin{theorem}\label{thm:VD}
Let $n>r\ge 1$, and let $\Delta=\VR(I_n;r)$. The following
conditions are equivalent:
\begin{enumerate}[label=\textup{(\roman*)}]
\item $\Delta$ is vertex decomposable;
\item $\Delta$ is Cohen--Macaulay over some field
(equivalently, over every field\textup);
\item $\Delta$ is pure;
\item $r=1$ or $n=r+1$.
\end{enumerate}
Moreover, over every field, with $S=\Bbbk[x_\alpha:\alpha\in\F^n]$,
\[
  2\le\depth(S/I_\Delta)\le 2^r
  ~\text{and}~
  \pd_S(S/I_\Delta)\ge 2^n-2^r,
\]
and equality holds in the upper bound for the depth if $r=1$ or
$n=r+1$.
\end{theorem}

\begin{proof}
\emph{\textup{(iv)}$\Rightarrow$\textup{(i)} and
\textup{(ii)}.} Suppose first that $r=1$. Then $\Delta=I_n$ is a
connected one-dimensional complex. Every connected graph with at
least one edge, regarded as a one-dimensional complex, is vertex
decomposable. Indeed, we may argue by induction on the number of
vertices. A single edge is a simplex. Otherwise, let $v$ be a leaf of
a spanning tree. Then the deletion $\X\setminus v$ is connected and
has at least two vertices, so none of its vertices is isolated, and
the link of $v$ is a nonempty zero-dimensional complex, which is
vertex decomposable. Hence no face of the link is a facet of the
deletion, and $v$ is a shedding vertex. Moreover, every connected
one-dimensional complex is Cohen--Macaulay over every field.
Suppose next that $n=r+1$. Then $\Delta$ is the boundary of a
$2^r$-dimensional cross-polytope, that is, the join of $2^r$ copies
of $S^0$. It is vertex decomposable because joins preserve vertex
decomposability~\cite{PB}, and it is Gorenstein, and hence
Cohen--Macaulay, over every field.

\emph{\textup{(i)}$\Rightarrow$\textup{(iv)}.} Suppose that
$r\ge 2$ and $n\ge r+2$. We show that $0$ is not a shedding vertex.
Since translations act transitively on the vertices of $\Delta$,
this shows that $\Delta$ has no shedding vertex; as $\Delta$ is not
a simplex, it is then not vertex decomposable.
The Hamming ball $B(0,1)=\{0,e_1,\ldots,e_n\}$ has diameter
$2\le r$. Choose a facet $\sigma$ of $\Delta$ containing it, and let
$\tau=\sigma\setminus\{0\}$, so that $\tau\in\lk(\Delta,0)$. Suppose
that $\tau\cup\{y\}\in\Delta\setminus 0$ for some $y\notin\sigma$.
Since $\sigma$ is a facet and $\Delta$ is flag, $y$ is not $r$-close
to $0$, that is, $\wt(y)\ge r+1$. If $i\notin\supp(y)$ for some $i$,
then $d_H(y,e_i)=\wt(y)+1>r$, contradicting $e_i\in\tau$. Hence
$y=\one$; but then $d_H(y,e_i)=n-1>r$ for every $i$, which is again a
contradiction. Thus $\tau$ is a facet of $\Delta\setminus 0$, and $0$
is not a shedding vertex.

\emph{\textup{(ii)}$\Rightarrow$\textup{(iii)}.} Every
Cohen--Macaulay complex is pure.

\emph{\textup{(iii)}$\Rightarrow$\textup{(iv)}.} Consider the
$r$-dimensional coordinate subcube
$Q_0=\{(\alpha,0^{n-r}):\alpha\in\F^r\}$. It is a facet of $\Delta$:
if $y=(\beta,\delta)\notin Q_0$, then $\delta\ne 0$ and
$d_H\bigl(y,(\overline\beta,0^{n-r})\bigr)=r+\wt(\delta)>r.$
Thus $\Delta$ has a facet with $2^r$ elements. Now suppose that
$r\ge 2$ and $n\ge r+2$, and write $r=2t$ or $r=2t+1$ with $t\ge 1$.
Kleitman's formula and the symmetry of binomial coefficients give
\[
  \kappa(r+2,r)=
  \begin{cases}
    2^{r+1}-\frac12\binom{2t+2}{t+1} & \text{if } r=2t,\\[2pt]
    2^{r+1}-\binom{2t+2}{t+1} & \text{if } r=2t+1,
  \end{cases}
\]
so the inequality $\kappa(r+2,r)>2^r$ is
equivalent to $\binom{2t+2}{t+1}<2^{2t+1}$. This holds because
$\binom{2k}{k}<2^{2k-1}$ for every $k\ge 2$: indeed,
$\binom42=6<8$, and
$\binom{2k+2}{k+1}/\binom{2k}{k}=2(2k+1)/(k+1)<4$. Since $\kappa$ is
increasing in $n$, we obtain $\kappa(n,r)\ge\kappa(r+2,r)>2^r$, so
$\Delta$ also has a facet with more than $2^r$ elements. Hence
$\Delta$ is not pure.

Finally, let $\Bbbk$ be any field. The $1$-skeleton of $\Delta$
contains the connected graph $I_n$, so $\Delta$ is connected and has
positive dimension; hence $\depth(S/I_\Delta)\ge 2$~\cite{BH}. On the
other hand, $Q_0$ is a facet, so $\lk(\Delta,Q_0)=\{\emptyset\}$ has
nonzero reduced homology in degree $-1$. By Hochster's formula for
local cohomology~\cite{BH}, this gives
$H^{|Q_0|}_{\mathfrak m}(S/I_\Delta)\ne 0$, and therefore
$\depth(S/I_\Delta)\le|Q_0|=2^r$. Since $S$ has $2^n$ variables, the
Auslander--Buchsbaum formula yields
$\pd_S(S/I_\Delta)=2^n-\depth(S/I_\Delta)\ge 2^n-2^r$. If $r=1$ or
$n=r+1$, then $\Delta$ is Cohen--Macaulay of dimension $2^r-1$, so
$\depth(S/I_\Delta)=\dim(S/I_\Delta)=2^r$.
\end{proof}

Theorem~\ref{thm:VD} shows that vertex decomposability cannot yield
new cases of Conjecture~\ref{conj:shukla}: when $r\ge 2$ and
$n\ge r+2$, no vertex of $\VR(I_n;r)$ is a shedding vertex.
Nevertheless, the equality $\theta(\VR(I_n;r))=2^r$ holds for
$n=r+2$ by Theorem~\ref{thm:VR-n-2}, which shows that the
$\theta$-number can be sharp even when vertex decomposability fails.

\subsection{Proof of Theorem~\texorpdfstring{\ref{thmC}}{C}}
\label{subsec:proofD}

\begin{proof}
Part~(i) is Theorem~\ref{thm:topology}(b), and part~(ii) is
Theorem~\ref{thm:topology}(c). Part~(iii) is
Proposition~\ref{prop:retraction}, and part~(iv) is
Theorem~\ref{thm:nonvanishing}(a). Part~(v) is
Theorem~\ref{thm:VD}.
\end{proof}

\vspace*{1mm}
\noindent\textbf{Acknowledgments.}
The author gained valuable insights into several topics related to this
paper through the workshop "Cohen--Macaulay Simplicial Complexes in
Graph Theory", held at the Chennai Mathematical Institute (CMI) in
2023, and the "Fourth Meru Combinatorics Conference", held at IIT Bhilai
in 2026. The author is grateful to the organizers of both events for
providing stimulating academic environments and opportunities for
learning and discussion. The author was partially supported by the
National Board for Higher Mathematics (NBHM) and the Science and
Engineering Research Board (SERB).

\vspace*{1mm}
\noindent\textbf{Data availability statement.}
Data sharing is not applicable to this article, as no datasets
were generated or analyzed during the current study.

\vspace*{1mm}
\noindent\textbf{Conflict of interest.}
The author declares that there are no competing interests
that could have influenced the work reported in this paper.

\vspace*{1mm}
\noindent\textbf{Declaration of AI-assisted writing.}
During the preparation of this manuscript, the author used
OpenAI's ChatGPT to improve the language, clarity, and
presentation of the text. The author subsequently reviewed
and edited the AI-assisted content and takes full
responsibility for the final version of the manuscript.

\bibliographystyle{abbrv}
\bibliography{ref}
\end{document}